\documentclass[11pt]{article}

\usepackage{amsmath,amsfonts,amssymb}
\usepackage{bbm}
\usepackage{cases}
\usepackage{color}
\usepackage[hidelinks]{hyperref}
\usepackage{amsthm}
\usepackage{changes}
\usepackage{pdfpages}

\def\|{\vert\vert}

\def\<{{\langle}}
\def\>{{\rangle}}

\theoremstyle{plain}
\newtheorem{theorem}{Theorem}[section]
\newtheorem{lemma}[theorem]{Lemma}
\newtheorem{corollary}[theorem]{Corollary}

\theoremstyle{definition}
\newtheorem{example}[theorem]{Example}
\newtheorem{definition}[theorem]{Definition}
\newtheorem{remark}[theorem]{Remark}

\allowdisplaybreaks

\begin{document}

\title{\vspace*{-1.5cm}
Well-posedness of stochastic partial differential equations with fully local monotone coefficients}

{\footnotesize
\author{Michael R\"{o}ckner$^{1}$,\ \ Shijie Shang$^{2}$,\ \ Tusheng Zhang$^{3}$\\
 {\em $^1$ Department of Mathematics, Bielefeld University,}\\
 {\em D-33501 Bielefeld, Germany. roeckner@math.uni-bielefeld.de}\\
 {\em $^2$ School of Mathematical Sciences,}\\
 {\em University of Science and Technology of China,}\\
 {\em Hefei, 230026, China. sjshang@ustc.edu.cn}\\
 {\em $^3$ School of Mathematics, University of Manchester,}\\
 {\em Oxford Road, Manchester, M13 9PL, UK. tusheng.zhang@manchester.ac.uk}
}}

\maketitle

\begin{abstract}
Consider stochastic partial differential equations (SPDEs)  with fully local monotone coefficients in a Gelfand triple $V\subseteq H \subseteq V^*$:
\begin{numcases}{}
 dX(t) = A(t,X(t))dt + B(t,X(t))dW(t),  \quad t\in (0,T], \nonumber\\
~ X(0) = x\in H, \nonumber
\end{numcases}
where
\begin{align*}
	A: [0,T]\times V \rightarrow V^* , \quad B: [0,T]\times V \rightarrow L_2(U,H)
\end{align*}
are measurable maps, $L_2(U,H)$ is the space of Hilbert-Schmidt operators from $U$ to $H$ and $W$ is a $U$-cylindrical Wiener process.
Such SPDEs include many interesting models in applied fields like fluid dynamics etc. In this paper, we establish
the well-posedness of the above SPDEs under fully local monotonicity condition solving a longstanding open problem. The conditions on the diffusion coefficient $B(t,\cdot)$ are allowed to depend on both the $H$-norm and $V$-norm. In the case of classical SPDEs, this means that $B(\cdot,\cdot)$ could also depend on the gradient of the solution.
The well-posedness is obtained through a combination of pseudo-monotonicity techniques and compactness arguments.

\end{abstract}


\medskip

\newpage
\tableofcontents

\section{Introduction}

\quad Let $H$ be a separable Hilbert space with inner product $(\cdot, \cdot)$ and norm $\|\cdot\|_{H}$. Let $V$ be a reflexive Banach space that is continuously and densely embedded into $H$. The norms of $V$ and its dual space $V^*$ are denoted by $\|\cdot\|_{V}$ and $\|\cdot\|_{V^*}$ respectively.
If we identify the Hilbert space $H$ with its dual space $H^*$ by the Riesz representation, then we obtain a Gelfand triple
\begin{align*}
V\subseteq H\subseteq V^*.
\end{align*}

\noindent We denote by $\langle f,v\rangle$ the dual pairing between $f\in V^*$ and $v\in V$. It is easy to see that
\begin{align}\label{eq P4 star}
(u,v)=\langle u,v\rangle,\ \ \ \forall\,u\in H,\ \ \forall\,v\in V.
\end{align}

\noindent Let $W$ be a cylindrical Wiener process on another separable Hilbert space $U$ defined on some probability space $(\Omega, \mathcal{F}, \mathbb{P})$ with normal filtration $\{\mathcal{F}_t\}$.


Let $T>0$ be fixed in this paper. Consider the following stochastic partial differential equations (SPDEs),
\begin{numcases}{}
 dX(t) = A(t,X(t))dt + B(t,X(t))dW(t),  \quad t\in (0,T], \nonumber\\
\label{1.a} ~ X(0)= x\in H,
\end{numcases}
where
\begin{align}
	A: [0,T]\times V \rightarrow V^* , \quad B: [0,T]\times V \rightarrow L_2(U,H)
\end{align}
are measurable maps. Here $L_2(U,H)$ is the space of Hilbert-Schmidt operators from $U$ to $H$ with the norm denoted by $\Vert \cdot\Vert_{L_2}$.

In this paper, we are concerned with the existence and uniqueness of variational solutions for the above SPDEs in the Gelfand triple. 

\vskip 0.3cm

In variational approach, monotone operators play a key role.
The theory of monotone operators was initiated in the fundamental work of Minty \cite{M1962}, and then studied systematically by Browder, and later developed by Leray and Lions, Hartman and Stampacchia, We refer the reader to \cite{Lions, Z1990} for a detailed exposition. 

\vskip 0.3cm


In the case of SPDEs, the monotone method was initially introduced by Pardoux in his pioneering work \cite{P1972, P1975}. Celebrated work was then carried out by Krylov and Rozovskii \cite{KR1979} and Gyongy \cite{G1982}.
We refer to \cite{PR, RRW, RW, RZ, Z} and references therein for the early applications of the variational approach to SPDEs. Now, there exists an enormous literature on the well-posedness of SPDEs within the variational framework, some of which we like to mention here.
Please see \cite{LR3, Liu2, EV, NS, GHV2022} for SPDEs with generalized coercivity conditions or Lyapunov conditions, \cite{BLZ, NTT, KR2021} for SPDEs driven by Levy noise,
\cite{{GRZ2009}} for martingale solutions of SPDEs,
 and \cite{BDR2016} for stochastic porous media equations as well as other  references.

\vskip 0.3cm

In 2010, the classical framework of the variational approach was substantially extended by Liu and R\"{o}ckner \cite{LR} for SPDEs with coefficients satisfying local monotonicity conditions, more precisely, for $u, v\in V$,
\begin{align}\label{220326.1753}
	\langle A(t,u) - A(t,v), u - v \rangle \leq [C+\rho(u)+\eta(v)]\Vert u - v \Vert_{H}^2 ,
\end{align}
where $\rho(u)$ or $\eta(v)$ are functions on the smaller space $V$, which are bounded on $V$-balls. However, in \cite{LR} it was required that only one of $\rho(u)$ and $\eta(v)$ is non-zero, namely, either  $\rho(u)\equiv 0$ or $\eta(v)\equiv 0$.
Nevertheless,  many interesting examples such as stochastic 2D Navier–Stokes equations can be included into this framework.
Later in 2011, Liu \cite{Liu} (see also \cite{LR2}) studied SPDEs satisfying a more general type of locally monotonicity conditions, i.e. both $\rho(u)$ and $\eta(v)$ in (\ref{220326.1753}) are not zero.
Liu used pseudo-monotone operators within the variational approach, but only SPDEs with additive noise were solved.
Since then, the well-posedness of SPDEs driven by multiplicative noise with fully  local monotone coefficients has been left as an open problem, which was mentioned a number of times in the  literature (see e.g. \cite{LR, Liu, LR2, LR3, BLZ, Liu2}).

\vskip 0.3cm

The purpose of this paper is to establish the well-posedness of SPDEs driven by multiplicative noise with fully  local monotone coefficients, solving this problem in the field, which has been open for some time. The notion of  pseudo-monotone operators plays an important role.
Pseudo-monotone operators (see Definition \ref{220405.2112} in Section \ref{220406.2023}) were initially introduced by Br\'{e}zis in \cite{B1968}, and further developed by many authors, see \cite{Lions, Z1990, VZ2019, VZ2021} and references therein.
Pseudo-monotone operators constitute a more general class than monotone operators.
In their works \cite{VZ2019, VZ2021}, Vallet and Zimmermann considered a class of nonlinear diffusion-convection equations with multiplicative and additive noise. It is worth noting that the nonlinear diffusion-convection operator considered there is not monotone/locally monotone but rather pseudo-monotone. They employed a semi-implicit Euler–Maruyama time discretization method to address the well-posedness of the solutions.

\vskip 0.3cm

Now we describe our approach in detail.
We distinguish two cases depending on  whether the diffusion coefficient $B$ is continuous on the Hilbert space $H$ or $B(t,u)$ is allowed to depend on the gradient $\nabla u$ of the solution function $u$. In Part I, we treat the case where $B$ is continuous on $H$. We will combine compactness arguments with techniques from the theory of pseudo-monotone operators.
Firstly, we establish some improved uniform moment estimates for the Galerkin approximating solutions $\{ X_n \}$ (see (\ref{220608.2019}) below).
Secondly, we prove the tightness of the laws  of $\{X_n\}$ in the space $C([0,T],V^*)\cap L^{\alpha}([0,T],H)$.
Thirdly, we apply the Prohorov theorem and the Skorokhod representation theorem to obtain that on some new probability space, $X_n$ converges almost surely to some element $X$ in $C([0,T],V^*)\cap L^{\alpha}([0,T],H)$ (along a subsequence).
Finally, we show that $B(\cdot, X_n(\cdot))$ converges strongly to $B(\cdot, X(\cdot))$ and that $A(\cdot, X_n(\cdot))$ converges weakly to $A(\cdot, X(\cdot))$. To this end, it is essential to show that $Y(\cdot)\rightarrow A(\cdot, Y(\cdot))$ is pseudo-monotone in certain $L^{\alpha}$-spaces.
Hence $X$ will be a probabilistically weak solution. The existence of probabilistically strong solutions follows from the pathwise uniqueness and the Yamada-Watanabe theorem.
In part II, we deal with the case where $B$ is allowed to depend on the gradient of the solution. We shall modify the  steps from Part I.
We  establish the tightness of the laws of $\{X_n\}$ in the space $L^{\alpha}([0,T],H)$. To identify every limit point $X$ as a solution of the SPDE (\ref{1.a}),
we use in an essential way the fact that the $X_n$ converges almost surely to some element $X$ in the space $L^{\alpha}([0,T],H)$ on some new probability space. Monotonicity techniques also play an important role. 

\vskip 0.3cm

The results of this paper can be applied to establish the existence and uniqueness of solutions for many interesting  stochastic nonlinear evolution equations. Examples are provided in Section \ref{220406.2046}. It should be pointed out that all the examples considered in \cite{PR, LR2, LR, Liu} can be covered by our framework, including the 2D Navier-Stokes equations, porous media equations, fast-diffusion equations, $p$-Laplacian equations, Burgers equations, Allen-Cahn equations, 3D Leray-$\alpha$ model, 2D Boussinesq system, 2D magneto-hydrodynamic equations, 2D Boussinesq model for the B{\'e}nard convection, 2D magnetic B{\'e}nard equations, some shell models of turbulence (GOY, Sabra, dyadic), power law fluids, the Ladyzhenskaya model, and the Kuramoto-Sivashinsky equations. Moreover, our main results are also applicable to the 3D tamed Navier-Stokes equations, some quasilinear PDEs, Cahn-Hilliard equations, liquid crystal models and Allen-Cahn-Navier-Stokes systems, which are not covered by the framework  in \cite{LR, LR2, PR}.

\vskip 0.3cm

Here are some conventions on constants. Throughout the paper, $C$ denotes a generic positive constant whose value may change from line to line. All other constants will be denoted by $C_1$, $C_2$,$\cdots$. They are all positive but their values are not important. The dependence of constants on parameters if needed will be indicated, e.g. $C_p, C_{\varepsilon}$.

\section{Part I }\label{220406.2023}
\setcounter{equation}{0}
\subsection{Hypotheses and Main Results}

Let us first recall the definition of pseudo-monotone operators.
\begin{definition}\label{220405.2112}
	An operator $A$ from $V$ to $V^*$ is said to be pseudo-monotone, if the following property holds: if $u_n$ converges weakly to $u$ in $V$ and
	\begin{align}
		\liminf_{n\rightarrow\infty} \langle A(u_n), u_n-u \rangle \geq 0 ,
	\end{align}
	then
	\begin{align}
		\limsup_{n\rightarrow\infty} \langle A(u_n), u_n -v \rangle \leq \langle A(u), u-v \rangle , \quad \forall\, v\in V .
	\end{align}
\end{definition}
\begin{remark}\label{220119.1500}
	If $A$ is a bounded operator from $V$ to $V^*$, i.e. $A$ maps every bounded set of $V$ to a bounded set of $V^*$, then the pseudo-monotonicity of $A$ is equivalent to the following property: if $u_n$ converges weakly to $u$ in $V$ and
	\begin{align}
		\liminf_{n\rightarrow\infty} \langle A(u_n), u_n-u \rangle \geq 0 ,
	\end{align}
	then $A(u_n)$ converges to $A(u)$ in the weak-* topology of $V^*$ and
	\begin{align}
		\lim_{n\rightarrow\infty} \langle A(u_n), u_n \rangle = \langle A(u), u \rangle .
	\end{align}
We refer the reader to Proposition 27.7 in \cite{Z1990} or Remark 5.2.12 in \cite{LR2}.
\end{remark}

\vskip 0.3cm

We introduce the following conditions on the coefficients $A$ and $B$.
Let $f\in L^1([0,T],\mathbb{R}_+)$ and $\alpha\in (1,\infty)$.
\begin{itemize}
	\item [(H1)] Hemicontinuity:  for a.e. $t\in [0,T]$, the map $\mathbb{R}\ni\lambda\longmapsto\langle A(t,u+\lambda v),x\rangle \in\mathbb{R}$ is continuous for any $u,v,x\in V$.
 	\item [(H2)] Local monotonicity: there exist nonnegative constants $\gamma$ and $C$ such that for a.e. $t\in [0,T]$ and any $u,v\in V$,
 \begin{gather}
 \begin{aligned}
& 2 \langle A(t,u)-A(t,v), u-v \rangle + \Vert B(t,u)-B(t,v)\Vert_{L_2}^2 \nonumber\\
\leq & [f(t)+ \rho(u)+\eta(v) ] \Vert u-v\Vert_{H}^2 , \\
\end{aligned}\\
\label{1227.2214} |\rho(u)|+|\eta(u)| \leq C (1+ \Vert u\Vert_V^{\alpha} )(1+\Vert u\Vert_H^{\gamma}),
\end{gather}
where $\rho$ and $\eta$ are two measurable functions from $V$ to $\mathbb{R}$.
	\item [(H2)$^\prime$] General local monotonicity: for any $R>0$, there exists a function $K_{R}\in L^1([0,T],\mathbb{R}_{+})$ such that for a.e. $t\in [0,T]$ and any $u,v\in V$ with $\Vert u\Vert_V \vee \Vert u\Vert_V \leq R$,
 \begin{align}
 \langle A(t,u)-A(t,v), u-v \rangle  \leq K_{R}(t) \Vert u-v\Vert_{H}^2  .
 \end{align}

 \end{itemize}
\begin{remark}
	Obviously, (H2)$^\prime$ is weaker than  (H2). It turns out that the assumption (H2)$^\prime$ is sufficient  for the existence of solutions, while (H2) is used for the pathwise uniqueness of solutions.
\end{remark}

 \begin{itemize}
	\item [(H3)] Coercivity: there exists a constant $c>0$ such that for a.e. $t\in [0,T]$ and any $u\in V$,
		\begin{align}
			2 \langle A(t,u),u \rangle +\Vert B(t,u)\Vert^2_{L_2} \leq f(t) (1+\Vert u\Vert_{H}^2) - c \Vert u \Vert_{V}^{\alpha} .
		\end{align}
	\item [(H4)] Growth: there exist nonnegative constants $\beta$ and $C$ such that for a.e. $t\in [0,T]$ and any $u\in V$,
		\begin{align}
			\Vert A(t,u)\Vert_{V^*}^{\frac{\alpha}{\alpha-1}} \leq \left(f(t)+ C \Vert u\Vert_V^{\alpha}\right)(1+\Vert u\Vert_H^{\beta}) .
		\end{align}
	\item [(H5)] For a.e. $t\in [0,T]$ and any sequence $\{u_n\}_{n=1}^{\infty}$ and $u$ in $V$ satisfying $\Vert u_n-u\Vert_{H}\rightarrow 0$,
\begin{align}\label{220124.1129}
	\Vert B(t,u_n)- B(t,u)\Vert_{L_2}\rightarrow 0 .
\end{align}
Moreover, there exists $g \in L^{1}([0,T],\mathbb{R}_+)$
such that for a.e. $t\in [0,T]$ and any $u\in V$,
\begin{align}\label{1230.1109}
	\Vert B(t,u)\Vert^2_{L_2} \leq g(t) (1+ \Vert u\Vert_{H}^2 ) .
\end{align}

\end{itemize}

\begin{remark}
In many applications, the coefficient $B$ is assumed to be locally Lipschitz and of linear growth with respect to $u$ in $H$-norm. So, (H5) is satisfied. The case where $\Vert B(t,u)\Vert_{L_2}$ also allows $\alpha$-growth in the $V$-norm will be studied in Section \ref{220123.1616}.
\end{remark}

Let us recall the following definition of variational solutions to stochastic partial differential equation (\ref{1.a}).

\begin{definition}\label{220123.1652}
An $H$-valued continuous and adapted stochastic process $(X(t))_{t\in [0,T]}$ is called a solution to equation (\ref{1.a}), if for its $dt \otimes \mathbb{P}$-equivalence class $\hat{X}$ we have
\begin{align*}
	\hat{X}\in L^{\alpha}([0,T]\times\Omega,dt\otimes \mathbb{P},V) \cap L^2([0,T]\times\Omega,dt\otimes \mathbb{P},H)
\end{align*}
with $\alpha$ in (H3) and $\mathbb{P}$-a.s.
\begin{align*}
	X(t)=X_0 + \int_0^t A(s,\bar{X}(s)) ds + \int_0^t B(s,\bar{X}(s)) dW(s) , \quad\forall\  t\in [0,T],
\end{align*}
where $\bar{X}$ is any $V$-valued progressively measurable $dt \otimes \mathbb{P}$-version of $\hat{X}$.
\end{definition}
%
%

Our main results in this section read as follows.
\begin{theorem}\label{1227.2229}
Suppose that the embedding $V\subseteq H$ is compact and (H1) (H2)$^{\prime}$ (H3) (H4) (H5) hold.
Then for any initial value $x\in H$, there exists a probabilistically weak solution to equation (\ref{1.a}). Furthermore, for any $p\geq 2$, the following moment estimate holds:
\begin{align}\label{1226.2116}
\mathbb{E} \Big\{\sup_{t\in [0,T]} \Vert X(t)\Vert_H^p\Big\} +  \mathbb{E}\left\{ \Big( \int_0^T \Vert X(s)\Vert_V^{\alpha} ds \Big)^{\frac{p}{2}}  \right\} < \infty.
\end{align}
Moreover, if (H2) is satisfied, then pathwise uniqueness holds for solutions of equation (\ref{1.a}) and hence there exists a unique probabilistically strong solution to equation (\ref{1.a}).

\end{theorem}

From the proof of Theorem \ref{1227.2229}, we obtain the following corollary.
\begin{corollary}\label{220310.1121}
Assume the embedding $V\subseteq H$ is compact, the operator $A(t,\cdot)$ is pseudo-monotone for a.e. $t\in [0,T]$,  and (H1) (H3) (H4) (H5) hold.
Then for any initial value $x\in H$, there exists a probabilistically weak solution to equation (\ref{1.a}), and estimate (\ref{1226.2116}) holds.
\end{corollary}

\begin{theorem}\label{220119.1651}

Suppose that the embedding $V\subseteq H$ is compact and (H1) (H2) (H3) (H4) (H5) hold.
Let $\{ x_n \}_{n=1}^{\infty}$ and $x$ be a sequence with $ \Vert x_n -x \Vert_H \rightarrow 0$.
Let $X(t,x)$ be the unique solution of (\ref{1.a}) with the initial value $x$.
Then  for any $p\geq 2$,
\begin{align}\label{1231.1417}
	 \lim_{n\rightarrow\infty}\mathbb{E}\Big[ \sup_{t\in [0,T]}\Vert X(t, x_n) - X(t, x) \Vert_H^p \Big] = 0 .
\end{align}

\end{theorem}

\begin{remark}
Compared with the local monotonicity condition used in \cite{LR}, the major difference is that in condition (H2) both $\rho$ and $\eta$ can be nonzero. In fact, in \cite{LR} it is required that $\rho(u)+\eta(v)$ either only depend on $u$ or only depend on $v$ when the equation is driven by multiplicative noise. This requirement was crucially used in \cite{LR}.
\end{remark}

\subsection{Proofs of the Main Results}\label{220120.2019}

In this section, we will prove Theorem \ref{1227.2229} and Theorem \ref{220119.1651}. Throughout this part, we will assume that the embedding $V\subseteq H$ is compact and (H1) (H2)$^{\prime}$ (H3) (H4) (H5) hold.


We will construct approximating solutions using the Galerkin method and then establish  the tightness of the laws of the approximating solutions in an appropriate space in order to obtain the existence of probabilistically weak solutions.

\vskip 0.6cm

Let $\{e_i\}_{i=1}^{\infty} \subset V$ be an orthonormal basis of $H$.
Let $H_n$ denote the $n$-dimensional subspace of $H$ spanned by $\{e_1,\dots, e_n\}$. Let $P_n: V^*\rightarrow H_n$ be defined by
\begin{align}\label{1229.2148}
  P_n u := \sum_{i=1}^{n}\langle u,e_i\rangle e_i .
\end{align}
Clearly, $P_n|_H$ is just the orthogonal projection of $H$ onto $H_n$.
Let $\{h_i\}_{i=1}^{\infty}$ be an orthonormal basis of Hilbert space $U$.
Set
\begin{align}\label{1229.2149}
	W_n(t) = Q_n W(t):= \sum_{i=1}^n \langle W(t), h_i \rangle h_i ,
\end{align}
where $Q_n$ is the orthogonal projection onto $\mathrm{span}\{h_1,\cdots,h_n\}$ in $U$.

For any integer $n\geq 1$, we consider the following stochastic differential equation in the finite-dimensional space  $H_n$,
\begin{align}\label{1228.2225}
	Y_n(t) = P_n x +\int_0^t P_n A(s,Y_n(s))ds + \int_0^t P_n B(s,Y_n(s))Q_n dW(s) .
\end{align}
Note that, under our considerations, the operator $u\in V\rightarrow A(t,u)\in V^*$ is demicontinuous for a.e. $t\in [0,T]$, that is for $u_n\rightarrow u$ in $V$, we have $A(t,u_n)\rightarrow A(t,u)$ weakly in $V^*$.
In fact, if the conditions (H1) and (H2)$^\prime$ hold, then we can apply Remark 4.1.1 in \cite{LR2} or Proposition 26.4 in \cite{Z1990}. If the condition $(H4)$ holds  and if  $A(t,\cdot)$ is pseudo-monotone for a.e. $t\in [0,T]$, then we can apply Proposition 27.7 in \cite{Z1990}.
Therefore, in both cases, the map $H_n \ni u \rightarrow P_n A(t,u) + P_n B(t,u)h(t) \in H_n$ is continuous for a.e. $t\in [0,T]$.
Combining this fact with the condition (H3), it's well-known that there exists a global probabilistically weak solution to the above stochastic differential equation (\ref{1228.2225}) (please refer to \cite{IW,LR2,HS}).
In the followings, we will establish uniform moment estimates for $\{ Y_n \}$. Thus, for simplicity of notations, we can think that  $\{Y_n\}$ still lives in the original probability space $(\Omega,\mathcal{F},\mathbb{P})$.


\begin{lemma}\label{1231.1011}
For any $p\geq 2$, there exists a constant $C_p$ such that
	\begin{align}\label{220608.2019}
		& \sup_{n\in\mathbb{N}} \bigg\{ \mathbb{E} \Big[ \sup_{t\in [0,T]} \Vert Y_n(t)\Vert_H^p \Big] + \mathbb{E} \Big( \int_0^T \Vert Y_n(t)\Vert_V^{\alpha} dt \Big)^{\frac{p}{2}}\bigg\} < C_p(1+\Vert x \Vert_H^p).
	\end{align}	
\end{lemma}
\begin{remark}
	The above improved estimates (compared to Lemma 2.2 in \cite{LR} and Lemma 4.2.9 in \cite{PR}) are crucial in the analysis below. In \cite{GHV2022}, the authors obtained a similar estimate as above under a different set of conditions.
\end{remark}
\noindent {\bf Proof}. It suffices to prove this lemma for large $p$. By Ito's formula,
\begin{align}\label{220117.1455}
	\Vert Y_n(t) \Vert_H^2 = & \Vert P_n x \Vert_H^2 + \int_0^t \Big[  2\langle A(s, Y_n(s)), Y_n(s) \rangle +\Vert P_n B(s,Y_n(s))Q_n \Vert_{L_2}^2 \Big] ds \nonumber\\
	& + 2\int_0^t ( Y_n(s), B(s,Y_n(s))Q_n dW(s) ),
\end{align}
here
$$\int_0^t ( Y_n(s), B(s,Y_n(s))Q_n dW(s) ):=\int_0^t ((B(s,Y_n(s))Q_n)^* Y_n(s), dW(s)).$$
Now, we apply Ito's formula to the real-valued process $\Vert Y_n(t) \Vert_H^2$ to get
\begin{align}\label{220118.1447}
	\Vert Y_n(t) \Vert_H^p = & \Vert P_n x \Vert_H^p + \frac{p}{2}\int_0^t \Vert Y_n(s) \Vert_H^{p-2} \Big[  2\langle A(s, Y_n(s)), Y_n(s) \rangle +\Vert P_nB(s,Y_n(s))Q_n \Vert_{L_2}^2 \Big] ds \nonumber\\
	& + \frac{p(p-2)}{2} \int_0^t \Vert Y_n(s) \Vert_H^{p-4}  \Vert (B(s,Y_n(s))Q_n)^* Y_n(s)\Vert_{U}^2  ds \nonumber\\
	& + p\int_0^t \Vert Y_n(s) \Vert_H^{p-2} ( Y_n(s), B(s,Y_n(s))Q_n dW(s) ) .
\end{align}
In view of (H3) and (H5), it follows from (\ref{220118.1447}) that there exist some positive constants $c, C$ such that
\begin{align}
	& \Vert Y_n(t) \Vert_H^p + \frac{pc}{2} \int_0^t \Vert Y_n(s) \Vert_{V}^{\alpha} \Vert Y_n(s) \Vert_{H}^{p-2} ds \nonumber\\
	& \leq \Vert P_n x \Vert_H^p + C \int_0^t [f(s) + g(s)] (1+\Vert Y_n(s)\Vert_H^2) \Vert Y_n(s) \Vert_H^{p-2} ds \nonumber\\
	& + p\int_0^t \Vert Y_n(s) \Vert_H^{p-2} ( Y_n(s), B(s,Y_n(s))Q_n dW(s) ) .
\end{align}
For $M>0$, we define the stopping time
\begin{align}
	\tau_{n,H}^M := \inf\{t\in [0,T]: \Vert Y_n(t)\Vert_H >M \}\wedge T
\end{align}
with the convention $\inf\emptyset = +\infty$. Then $\tau_{n,H}^{M}\rightarrow T$, $\mathbb{P}$-a.s. as $M\rightarrow\infty$, for every $n$.
Taking the supremum over $t\leq r\wedge\tau_{n,H}^M$ for any fixed $r\in [0,T]$ and then taking expectations on both sides of the above inequality yield
\begin{align}\label{220117.1434}
	& \mathbb{E} \Big[ \sup_{t\leq r\wedge\tau_{n,H}^M}\Vert Y_n(t) \Vert_H^p  \Big] + \frac{pc}{2} \mathbb{E}\int_0^{r\wedge\tau_{n,H}^M} \Vert Y_n(s) \Vert_{V}^{\alpha} \Vert Y_n(s) \Vert_{H}^{p-2} ds \nonumber\\
	 \leq &  \Vert x \Vert_{H}^p + C\int_0^T [f(s)+g(s)]ds + C\mathbb{E} \int_0^{r\wedge\tau_{n,H}^M} [f(s)+g(s)]  \Vert Y_n(s) \Vert_H^p ds \nonumber\\
	 & +  p \mathbb{E} \bigg\{ \sup_{t\leq r\wedge\tau_{n,H}^M}\bigg| \int_0^t \Vert Y_n(s) \Vert_H^{p-2} ( Y_n(s), B(s,Y_n(s))Q_n dW(s) ) \bigg|\bigg\} .
\end{align}
By the BDG, H\"older and Young inequalities, using (H5) we obtain for any $\varepsilon>0$,
\begin{align}\label{220117.1435}
	& p \mathbb{E} \bigg\{ \sup_{t\leq r\wedge\tau_{n,H}^M}\bigg| \int_0^t \Vert Y_n(s) \Vert_H^{p-2} ( Y_n(s), B(s,Y_n(s))Q_n dW(s) ) \bigg|\bigg\} \nonumber\\
	\leq &  C \mathbb{E} \bigg( \int_0^{r\wedge\tau_{n,H}^M} \Vert Y_n(s) \Vert_H^{2p-2} \Vert B(s,Y_n(s)) \Vert_{L_2}^2 ds \bigg)^{\frac{1}{2}} \nonumber\\
	\leq & C\mathbb{E} \bigg( \sup_{t\leq r\wedge\tau_{n,H}^M}\Vert Y_n(t)\Vert_H^p \cdot\int_0^{r\wedge\tau_{n,H}^M} \Vert Y_n(s) \Vert_H^{p-2} \Vert B(s,Y_n(s))  \Vert_{L_2}^2 ds \bigg)^{\frac{1}{2}} \nonumber\\
	\leq & \varepsilon \mathbb{E} \Big[ \sup_{t\leq r\wedge\tau_{n,H}^M}\Vert Y_n(t)\Vert_H^p \Big] + C_{\varepsilon} \mathbb{E} \int_0^{r\wedge\tau_{n,H}^M} \Vert Y_n(s) \Vert_H^{p-2} \Vert B(s,Y_n(s)) \Vert_{L_2}^2 ds \nonumber\\
	\leq & \varepsilon \mathbb{E} \Big[ \sup_{t\leq r\wedge\tau_{n,H}^M}\Vert Y_n(t)\Vert_H^p \Big] + C_{\varepsilon} \int_0^T g(s) ds + C_{\varepsilon} \mathbb{E}\int_0^{r\wedge\tau_{n,H}^M} g(s)\Vert Y_n(s) \Vert_H^{p} ds .
\end{align}
Combining (\ref{220117.1434}) and (\ref{220117.1435}) together, appropriately choosing the parameter $\varepsilon$ and applying Gronwall's inequality give
\begin{align}
	& \mathbb{E} \Big[ \sup_{t\leq r\wedge\tau_{n,H}^M}\Vert Y_n(t) \Vert_H^p  \Big] + C \mathbb{E}\int_0^{r\wedge\tau_{n,H}^M} \Vert Y_n(s) \Vert_{V}^{\alpha} \Vert Y_n(s) \Vert_{H}^{p-2} ds \nonumber\\
	 \leq &  C \bigg( \Vert x \Vert_{H}^p + \int_0^T [f(s)+g(s)]ds \bigg) \exp\bigg( C\int_0^T [f(s)+g(s)]ds \bigg) .
\end{align}
Letting $M\rightarrow\infty$ and applying Fatou's lemma, we obtain
\begin{align}\label{220117.1519}
		& \sup_{n\in\mathbb{N}}\left\{ \mathbb{E} \Big[ \sup_{t\in [0,T]}\Vert Y_n(t) \Vert_H^p  \Big] +  \mathbb{E}\int_0^{T} \Vert Y_n(s) \Vert_{V}^{\alpha} \Vert Y_n(s) \Vert_{H}^{p-2} ds \right\} < \infty.
\end{align}
Using (H3), it follows from (\ref{220117.1455}) that
\begin{align*}
		& \Vert Y_n(t) \Vert_H^2 + c \int_0^t \Vert Y_n(s)\Vert_V^{\alpha} ds \nonumber\\
		\leq & \Vert P_n x \Vert_H^2 + \int_0^t [f(s)+g(s)](1+\Vert Y_n (s)\Vert_H^2) ds \nonumber\\
	& + 2\int_0^t ( Y_n(s), B(s,Y_n(s))Q_n dW(s) ) .
\end{align*}
Hence
\begin{align}\label{220117.1520}
	\mathbb{E}\bigg(\int_0^t \Vert Y_n(s)\Vert_V^{\alpha} ds \bigg)^{\frac{p}{2}}
	\leq &  C\Vert x\Vert_H^{p} + C \mathbb{E}\bigg( \int_0^t [f(s)+g(s)](1+\Vert Y_n (s)\Vert_H^2) ds \bigg)^{\frac{p}{2}} \nonumber\\
	& + C\mathbb{E}\bigg| \int_0^t ( Y_n(s), B(s,Y_n(s))Q_n dW(s) ) \bigg|^{\frac{p}{2}} \nonumber\\
	\leq & C\Vert x\Vert_H^{p} + C \Big( 1+ \mathbb{E} \Big[ \sup_{s\in [0,T]}\Vert Y_n (s)\Vert_H^p \Big] \Big) \nonumber\\
	& + C\mathbb{E}\bigg| \int_0^t ( Y_n(s), B(s,Y_n(s))Q_n dW(s) ) \bigg|^{\frac{p}{2}} .
\end{align}
Again by the BDG inequality and (H5), we have
\begin{align}\label{220117.1521}
	& C\mathbb{E}\bigg| \int_0^t ( Y_n(s), B(s,Y_n(s))Q_n dW(s) ) \bigg|^{\frac{p}{2}} \nonumber\\
	\leq & C \mathbb{E} \bigg( \int_0^t \Vert Y_n(s) \Vert_H^2 \Vert B(s, Y_n(s))\Vert_{L_2}^2 ds \bigg)^{\frac{p}{4}} \nonumber\\
	\leq & C \mathbb{E} \bigg( \int_0^T g(s) ds \cdot\Big[ 1 +\sup_{s\in [0,T]} \Vert Y_n(s) \Vert_H^4 \Big]  \bigg)^{\frac{p}{4}} \nonumber\\
	\leq & C \Big( 1+ \mathbb{E} \Big[ \sup_{s\in [0,T]} \Vert Y_n(s) \Vert_H^p \Big] \Big) .
\end{align}
Combining (\ref{220117.1519})-(\ref{220117.1521}) together yields
\begin{align}
	\sup_{n\in\mathbb{N}} \mathbb{E}\bigg(\int_0^T \Vert Y_n(s)\Vert_V^{\alpha} ds \bigg)^{\frac{p}{2}}  < \infty .
\end{align}
The proof is complete.
$\blacksquare$

\vskip 0.3cm


For $M>0$, we define the stopping time
\begin{align}
	\tau_{n}^M := & \inf \{t\in [0,T]: \Vert Y_n(t)\Vert_H^2 > M \} \nonumber\\
  & \wedge \inf\Big\{t\in [0,T]: \int_0^t\Vert Y_n(s)\Vert_{V}^{\alpha}ds >M\Big\} \wedge T
\end{align}
with the convention $\inf\emptyset = +\infty$.
By the Chebyshev inequality,
Lemma \ref{1231.1011} implies that
\begin{align}\label{1229.1652}
	\lim_{M\rightarrow\infty}\sup_{n\in \mathbb{N}} \mathbb{P}(\tau_{n}^M <T) =0.
\end{align}

The next result gives the tightness of the laws of $\{Y_n\}$.
\vskip 0.3cm
\begin{lemma}\label{1229.2041}
	$\{Y_n\}_{n=1}^{\infty}$ is tight in the space $C([0,T],V^*) \cap L^{\alpha}([0,T],H)$.
\end{lemma}
\noindent {\bf Proof}.
It suffices to prove that $\{Y_n\}_{n=1}^{\infty}$ is tight in $C([0,T],V^*)$ and in $L^{\alpha}([0,T],H)$ separately.

We first show that $\{Y_n\}_{n=1}^{\infty}$ is tight in $C([0,T],V^*)$.
According to the definition of the stopping time $\tau_n^M$ and (\ref{1229.1652}), we have
\begin{align}\label{220117.1625}
	\lim_{M\rightarrow\infty} \sup_{n\in\mathbb{N}}\mathbb{P}\Big(\sup_{t\in [0,T]}\Vert Y_n(t)\Vert_H >\sqrt{M}\Big) \leq \lim_{M\rightarrow\infty} \sup_{n\in\mathbb{N}}\mathbb{P}(\tau_{n}^M <T) =0 .
\end{align}
Since the embedding $H \subseteq V^*$ is compact, using (\ref{220117.1625}) and the fact that the Skorokhod topology in $D([0,T],V^*)$ restricted to $C([0,T],V^*)$ coincides with the uniform topology given by the supremum norm,  by Theorem 3.1 in \cite{J} it is sufficient to show that for every $e\in P_m H$, $m\in\mathbb{N}$, $\{\langle Y_n, e\rangle \}_{n=1}^{\infty}$ is tight in the space $C([0,T],\mathbb{R})$.
By (\ref{220117.1625}) and the Aldou tightness criterion (see Theorem 1 in \cite{A}), it suffices to show that for any stopping time $0\leq \zeta_n\leq T$ and for any $\varepsilon>0$,
\begin{align}\label{1229.1657}
	\lim_{\delta\rightarrow 0}\sup_{n\in\mathbb{N}}\mathbb{P}\Big( |\langle Y_n(\zeta_n +\delta)-Y_n(\zeta_n) , e\rangle | >\varepsilon\Big) =0 ,
\end{align}
 where $\zeta_n +\delta := T\wedge(\zeta_n +\delta)\vee 0$.
Set $Y_n^M(t):= Y_n(t\wedge\tau_n^M)$. By the Chebyshev inequality, we have
\begin{align}\label{1229.1648}
	& \mathbb{P}\Big(| \langle Y_n(\zeta_n +\delta)-Y_n(\zeta_n) , e\rangle |>\varepsilon\Big) \nonumber\\
	\leq & \mathbb{P}\Big(| \langle Y_n(\zeta_n +\delta)-Y_n(\zeta_n) , e\rangle | >\varepsilon, \tau_n^M\geq T\Big) + \mathbb{P}(\tau_n^M<T) \nonumber\\
	\leq &  \frac{1}{\varepsilon^{\alpha}} \mathbb{E}| \langle Y_n^M(\zeta_n +\delta)-Y_n^M(\zeta_n) , e\rangle |^{\alpha} + \mathbb{P}(\tau_n^M<T).
\end{align}
By the equation (\ref{1228.2225}) and the BDG inequality, it follows that
\begin{align}\label{1229.1649}
	& \mathbb{E} | \langle Y_n^M(\zeta_n +\delta)-Y_n^M(\zeta_n) , e\rangle |^{\alpha} \nonumber\\
	\leq & 2^{\alpha-1} \,\mathbb{E}\left(\int_{\zeta_n\wedge\tau_n^M}^{(\zeta_n+\delta)\wedge\tau_n^M} | \langle P_n A(s,Y_n(s)) , e\rangle | ds\right)^{\alpha} \nonumber\\
	& + 2^{\alpha-1} \,\mathbb{E}\left(\int_{\zeta_n\wedge\tau_n^M}^{(\zeta_n+\delta)\wedge\tau_n^M}  \Vert e \Vert_{H}^2 \Vert P_n B(s,Y_n(s))Q_n \Vert_{L_2}^{2} ds\right)^{\frac{\alpha}{2}} \nonumber\\
	=& : I_n + II_n .
\end{align}
Since $e\in P_m H$, we have
\begin{align*}
	\sup_{n\in\mathbb{N}} \Vert P_n e\Vert_V <\infty.
\end{align*}
By H\"{o}lder's inequality, (H4) and the above inequality, it follows that
\begin{align}\label{1229.1650}
	I_n \leq &  C\,\mathbb{E}\left\{|\delta| \times\Big[\int_{\zeta_n\wedge\tau_n^M}^{(\zeta_n+\delta)\wedge\tau_n^M} | \langle A(s,Y_n(s)) , P_n e \rangle |^{\frac{\alpha}{\alpha-1}}ds\Big]^{\alpha-1}\right\} \nonumber\\
	\leq &  C |\delta| \, \mathbb{E}\left[\int_0^{T\wedge\tau_n^M} \Vert P_n e\Vert_V^{\frac{\alpha}{\alpha-1}}\big(f(s)+C\Vert Y_n(s)\Vert_V^{\alpha}\big)\big(1+\Vert Y_n(s)\Vert^{\beta}_H\big)ds\right]^{\alpha-1} \nonumber\\
	\leq & C_M |\delta|  .
\end{align}
Similarly, by (H5) we have
\begin{align*}
	II_n \leq & C \,\mathbb{E}\left(\int_{\zeta_n\wedge\tau_n^M}^{(\zeta_n+\delta)\wedge\tau_n^M} \Vert e\Vert_H^2 \, g(s) (1+ \Vert  Y_n(s)\Vert_{H}^2 ) ds\right)^{\frac{\alpha}{2}} \nonumber\\
	\leq & C_M  \,\mathbb{E}\left(\int_{\zeta_n\wedge\tau_n^M}^{(\zeta_n+\delta)\wedge\tau_n^M} g(s)  ds\right)^{\frac{\alpha}{2}} .
\end{align*}
Note that $g\in L^1([0,T],\mathbb{R}_+)$.  By the absolute continuity of the Lebesgue integral and the dominated convergence theorem, we get
\begin{align}\label{1229.1651}
	\lim_{\delta\rightarrow 0}\sup_{n\in\mathbb{N}} II_n = 0 .
\end{align}
%
%
Combining (\ref{1229.1649})-(\ref{1229.1651}) together yields
\begin{align}\label{1229.1737}
	 \lim_{\delta\rightarrow 0}\sup_{n\in\mathbb{N}} \mathbb{E} | \langle Y_n^M(\zeta_n +\delta)-Y_n^M(\zeta_n) , e\rangle |^{\alpha} = 0 . 	
\end{align}
In view of (\ref{1229.1652}) and (\ref{1229.1737}), letting $\delta\rightarrow 0$ and then $M\rightarrow \infty$ in (\ref{1229.1648}) yield (\ref{1229.1657}). Thus we complete the proof of the tightness of $\{Y_n\}_{n=1}^{\infty}$ in $C([0,T],V^*)$.
%
%
%

\vskip 0.3cm

Next, we prove that $\{Y_n\}_{n=1}^{\infty}$ is tight in $L^{\alpha}([0,T],H)$.
Since, by Lemma \ref{1231.1011}
\begin{align}
	\sup_{n\in\mathbb{N}}\mathbb{E} \int_0^T \Vert Y_n(t) \Vert_V^{\alpha} dt < \infty,
\end{align}
by Lemma \ref{220117.1949} in Appendix, it is sufficient to show that for any $\epsilon >0$,
\begin{align}\label{220117.2120}
		\lim_{\delta\rightarrow 0+} \sup_{n\in\mathbb{N}} \mathbb{P} \left( \int_0^{T-\delta} \Vert Y_n (t+\delta) - Y_n (t) \Vert_{H}^{\alpha} dt > \epsilon \right)  = 0 .
\end{align}
Set $Y_n^M(t):= Y_n(t\wedge\tau_n^M)$ as before and note that
\begin{align}\label{220123.1745}
	& \mathbb{P} \left( \int_0^{T-\delta} \Vert Y_n (t+\delta) - Y_n (t) \Vert_{H}^{\alpha} dt > \epsilon \right) \nonumber\\
	\leq & \mathbb{P} \left( \int_0^{T-\delta} \Vert Y_n (t+\delta) - Y_n (t) \Vert_{H}^{\alpha} dt > \epsilon , \tau_n^M \geq T\right) + \mathbb{P} \left( \tau_n^M <T \right) \nonumber\\
	\leq & \frac{1}{\epsilon} \mathbb{E} \int_0^{T-\delta} \Vert Y_n^M (t+\delta) - Y_n^M (t) \Vert_{H}^{\alpha} dt + \mathbb{P} \left( \tau_n^M <T \right) .
\end{align}
If we have proved that for any fixed $M>0$,
\begin{align}\label{220117.2124}
	\lim_{\delta\rightarrow 0+} \sup_{n\in\mathbb{N}} \mathbb{E} \int_0^{T-\delta} \Vert Y_n^M (t+\delta) - Y_n^M (t) \Vert_{H}^{\alpha} dt = 0,
\end{align}
then, in light of (\ref{1229.1652}), letting $\delta\rightarrow 0$ and then $M\rightarrow \infty$ in (\ref{220123.1745}) we have (\ref{220117.2120}), completing the proof of the tightness of $\{Y_n\}_{n=1}^{\infty}$ in $L^{\alpha}([0,T],H)$. Therefore, it remains  to prove (\ref{220117.2124}). To this end, we consider two cases according to the value of $\alpha$.

We first consider the case $1<\alpha\leq 2$. Using Ito's formula to the process $\{ Y_n(r) - Y_n(t\wedge\tau_n^M), r\in [t\wedge\tau_n^M,(t+\delta)\wedge\tau_n^M] \}$, we have
\begin{align}
	& \mathbb{E}\Vert Y_n^M (t+\delta) - Y_n^M (t) \Vert_{H}^2 \nonumber\\
	= & \mathbb{E}\int_{t\wedge\tau_n^M}^{(t+\delta)\wedge\tau_n^M} 2\langle A(r,Y_n(r)), Y_n(r) - Y_n(t\wedge\tau_n^M) \rangle dr \nonumber\\
	& + \mathbb{E}\int_{t\wedge\tau_n^M}^{(t+\delta)\wedge\tau_n^M} \Vert P_n B(r,Y_n(r))Q_n \Vert_{L_2}^2 dr .
\end{align}
It follows that
\begin{align}\label{220117.2220}
	& \mathbb{E} \int_0^{T-\delta} \Vert Y_n^M (t+\delta) - Y_n^M (t) \Vert_{H}^{2} dt \nonumber\\
	= & \mathbb{E} \int_0^{T-\delta} \int_{t\wedge\tau_n^M}^{(t+\delta)\wedge\tau_n^M} \big[ 2\langle A(r,Y_n(r)), Y_n(r)\rangle + \Vert P_n B(r,Y_n(r))Q_n \Vert_{L_2}^2 \big] drdt   \nonumber\\
	& - 2\mathbb{E}\int_0^{T-\delta}  \int_{t\wedge\tau_n^M}^{(t+\delta)\wedge\tau_n^M} \langle A(r,Y_n(r)), Y_n(t\wedge\tau_n^M) \rangle drdt \nonumber\\
	=& : I_1 + I_2 .
\end{align}
By the Fubini theorem and (H3) it follows that
\begin{align}\label{220118.1648}
	I_1 = &  \mathbb{E} \int_0^{T\wedge\tau_n^M} \big[ 2\langle A(r,Y_n(r)), Y_n(r) \rangle + \Vert P_n B(r,Y_n(r))Q_n \Vert_{L_2}^2 \big]  \int_{0\vee(r-\delta)}^{r} \mathbf{1}_{\{\tau_n^M>t\}}  dtdr \nonumber\\
	\leq & \delta \,  \mathbb{E} \int_0^{T\wedge\tau_n^M} f(s) (1+\Vert Y_n(s)\Vert_H^2) ds \nonumber\\
	\leq &  C_M \delta .
\end{align}
Applying the Fubini theorem and (H4) it follows that
\begin{align}\label{220117.2221}
	|I_2| = & 2 \bigg|\mathbb{E} \int_0^{T\wedge\tau_n^M}  \int_{0\vee(r-\delta)}^{r} \mathbf{1}_{\{\tau_n^M>t\}} \langle A(r,Y_n(r)), Y_n(t\wedge\tau_n^M) \rangle dt dr\bigg| \nonumber\\
	\leq &  2\mathbb{E} \int_0^{T\wedge\tau_n^M} \Vert A(r,Y_n(r)) \Vert_{V^*}  \int_{0\vee(r-\delta)}^{r} \Vert Y_n(t\wedge\tau_n^M) \Vert_V  dt dr \nonumber\\
	\leq & 2\delta^{\frac{\alpha-1}{\alpha}}  \mathbb{E} \int_0^{T\wedge\tau_n^M} \Vert A(r,Y_n(r)) \Vert_{V^*} dr \bigg(\int_0^{T\wedge\tau_n^M} \Vert Y_n(t) \Vert_V^{\alpha} dt \bigg)^{\frac{1}{\alpha}} \nonumber\\
	\leq & C_M \delta^{\frac{\alpha-1}{\alpha}} .
\end{align}
Combining (\ref{220117.2220})-(\ref{220117.2221}) together yields
\begin{align}
	\sup_{n\in\mathbb{N}} \mathbb{E} \int_0^{T-\delta} \Vert Y_n^M (t+\delta) - Y_n^M (t) \Vert_{H}^{2} dt \leq C_M (\delta + \delta^{\frac{\alpha-1}{\alpha}}).
\end{align}
By H\"{o}lder's inequality, we see that
\begin{align}
		& \lim_{\delta\rightarrow 0+} \sup_{n\in\mathbb{N}} \mathbb{E} \int_0^{T-\delta} \Vert Y_n^M (t+\delta) - Y_n^M (t) \Vert_{H}^{\alpha} ds \nonumber\\
		\leq & C \bigg( \lim_{\delta\rightarrow 0+} \sup_{n\in\mathbb{N}} \mathbb{E} \int_0^{T-\delta} \Vert Y_n^M (t+\delta) - Y_n^M (t) \Vert_{H}^{2} ds \bigg)^{\frac{\alpha}{2}} = 0.
\end{align}
Thus we have proved (\ref{220117.2124}) for the case: $\alpha\leq 2$.

Now, we consider the remaining case: $\alpha> 2$. Applying Ito's formula to the function $\Vert \cdot \Vert_H^{\alpha}$ and then taking expectations, we have
\begin{align}
		& \mathbb{E}\Vert Y_n^M (t+\delta) - Y_n^M (t) \Vert_H^{\alpha} \nonumber\\
		= & \frac{\alpha}{2}\mathbb{E}\int_{t\wedge\tau_n^M}^{(t+\delta)\wedge\tau_n^M} \Vert Y_n(r) - Y_n(t\wedge\tau_n^M) \Vert_H^{\alpha-2} \Big[ 2 \langle A(r,Y_n(r)), \nonumber\\
		&~~~~~~~~~Y_n(r) - Y_n(t\wedge\tau_n^M) \rangle +\Vert P_n B(r,Y_n(r))Q_n \Vert_{L_2}^2 \Big] dr \nonumber\\
	& + \frac{\alpha(\alpha-2)}{2} \mathbb{E}\int_{t\wedge\tau_n^M}^{(t+\delta)\wedge\tau_n^M}  \Vert Y_n(r) - Y_n(t\wedge\tau_n^M) \Vert_H^{\alpha-4} \nonumber\\
	& ~~~~~~~~~\times\Vert \big(B(r,Y_n(r))Q_n\big)^*\big(Y_n(r) - Y_n(t\wedge\tau_n^M)\big) \Vert_{U}^2 dr .
\end{align}
By the Fubini theorem, it follows that
\begin{align}\label{220118.1623}
	& \mathbb{E} \int_0^{T-\delta} \Vert Y_n^M (t+\delta) - Y_n^M (t) \Vert_{H}^{\alpha} dt \nonumber\\
	= & \frac{\alpha}{2}\mathbb{E} \int_0^{T-\delta} \int_{t\wedge\tau_n^M}^{(t+\delta)\wedge\tau_n^M} \Vert Y_n(r) - Y_n(t\wedge\tau_n^M) \Vert_H^{\alpha-2} \nonumber\\
	& ~~~~~~~~\times\Big[ 2\langle A(r,Y_n(r)), Y_n(r) \rangle + \Vert P_n B(r,Y_n(r))Q_n \Vert_{L_2}^2 \Big] dr dt \nonumber\\
	& - \alpha \mathbb{E} \int_0^{T-\delta} \int_{t\wedge\tau_n^M}^{(t+\delta)\wedge\tau_n^M} \Vert Y_n(r) - Y_n(t\wedge\tau_n^M) \Vert_H^{\alpha-2} \langle A(r,Y_n(r)), Y_n(t\wedge\tau_n^M) \rangle dr dt	\nonumber\\
	& + \frac{\alpha(\alpha-2)}{2} \mathbb{E}\int_0^{T-\delta} \int_{t\wedge\tau_n^M}^{(t+\delta)\wedge\tau_n^M}  \Vert Y_n(r) - Y_n(t\wedge\tau_n^M) \Vert_H^{\alpha-4} \nonumber\\
	& ~~~~~~~~~\times\Vert \big(B(r,Y_n(r))Q_n\big)^* \big(Y_n(r) - Y_n(t\wedge\tau_n^M)\big) \Vert_{U}^2 dr dt \nonumber\\
	=& : J_1 + J_2 + J_3 .
\end{align}
Similarly to (\ref{220118.1648}) and (\ref{220117.2221}), one can show that
\begin{align}
	J_1 \leq & C_M \delta , \\
	|J_2| \leq &  C_M \delta^{\frac{\alpha-1}{\alpha}} .
\end{align}
On the other hand, by (H5) and the Fubini theorem it follows that
\begin{align}\label{220118.1656}
	J_3 \leq & C_M \mathbb{E} \int_0^{T\wedge\tau_n^M} \Vert P_n B(r,Y_n(r))Q_n \Vert_{L_2}^2  \int_{0\vee(r-\delta)}^{r} \mathbf{1}_{\{\tau_n^M>t\}}  dt dr\nonumber\\
	\leq & C_M \delta .
\end{align}
Combining (\ref{220118.1623})-(\ref{220118.1656}) together, we obtain
\begin{align}
	\sup_{n\in\mathbb{N}} \mathbb{E} \int_0^{T-\delta} \Vert Y_n^M (t+\delta) - Y_n^M (t) \Vert_{H}^{\alpha} dt \leq C_M (\delta + \delta^{\frac{\alpha-1}{\alpha}}) .
\end{align}
Therefore, (\ref{220117.2124}) is proved. Thus, we complete the proof of this lemma.
$\blacksquare$

\vskip 0.6cm

Set
\[\Upsilon:= [ L^{\alpha}([0,T],H) \cap C([0,T],V^*)] \times C([0,T],U_1) ,\]
where $U_1$ is a Hilbert space such that the imbedding $U\subset U_1$ is Hilbert-Schmidt.
From Lemma \ref{1229.2041}, we know that the family of the laws $\mathcal{L}(Y_n, W)$ of the random vectors $(Y_n, W)$ is tight in $\Upsilon$. By the Prohorov theorem and the modified Skorokhod representation theorem (see Theorem A.1 in \cite{NTT} or Theorem C.1 in \cite{BHR}), there exist a new probability space $(\widetilde{\Omega}, \widetilde{\mathcal{F}},\widetilde{\mathbb{P}})$ and a sequence of $\Upsilon$-valued random vectors $\{(\widetilde{X}_n, \widetilde{W}_n)\}$ and $(\widetilde{X}, \widetilde{W})$ such that
\begin{itemize}
	\item [(i)] $\widetilde{W}_n=\widetilde{W}$ for any $n\in\mathbb{N}$, \ $\widetilde{\mathbb{P}}$-a.s. ;
	\item [(ii)] $\mathcal{L}(\widetilde{X}_n, \widetilde{W}_n)=\mathcal{L}(Y_n, W)$ ;
	\item [(iii)] $\widetilde{\mathbb{P}}$-a.s.,
		\begin{align}\label{220119.1617}
			\Vert \widetilde{X}_n -\widetilde{X}\Vert_{L^{\alpha} ([0,T],H)}+ \Vert \widetilde{X}_n -\widetilde{X}\Vert_{C([0,T], V^*)}\rightarrow 0  .
		\end{align}
\end{itemize}
Next, we will show that $(\widetilde{X}, \widetilde{W})$ is a solution to equation (\ref{1.a}).

\vskip 0.3cm

Let $\widetilde{\mathcal{F}}_t$ be the filtration satisfying the usual conditions and generated by
\begin{align*}
	\{\widetilde{X}_n(s), \widetilde{X}(s), \widetilde{W}(s) : s\in [0,t] \} .
\end{align*}
Then $\widetilde{W}$ is an $\{\widetilde{\mathcal{F}}_t\}$-cylindrical Wiener process on $U$.
From the equation (\ref{1228.2225}) satisfied by the random vector $(Y_n, W)$, it follows that
\begin{align}\label{1230.1649}
  \widetilde{X}_n(t)=& P_n x + \int_0^t P_n A(s,\widetilde{X}_n(s)) ds  +\int_0^t P_n B(s,\widetilde{X}_n(s))Q_n d\widetilde{W}(s), \quad t\in[0,T] .
\end{align}
Moreover, $\{\widetilde{X}_n\}$ also satisfies the same moment estimates as $\{Y_n\}$ in Lemma \ref{1231.1011}, i.e. for any $p\geq 2$,
\begin{align}\label{1230.1048}
	\sup_{n\in\mathbb{N}} \bigg\{ \widetilde{\mathbb{E}} \Big[ \sup_{t\in [0,T]} \Vert \widetilde{X}_n(t)\Vert_H^p \Big] + \widetilde{\mathbb{E}} \Big( \int_0^T \Vert \widetilde{X}_n(t)\Vert_V^{\alpha} dt \Big)^{\frac{p}{2}} \bigg\} <\infty.
\end{align}	
Since $\Vert \cdot \Vert_H$ and $\Vert \cdot \Vert_V$ are lower semicontinuous in $V^*$,
by (\ref{220119.1617}) and Fatou's lemma, we obtain
\begin{align}\label{1230.1117}
	  \widetilde{\mathbb{E}}\sup_{t\in[0,T]}\Vert  \widetilde{X}(t)\Vert^p_H \leq & \widetilde{\mathbb{E}}\sup_{t\in[0,T]}\liminf_{n\rightarrow\infty} \Vert \widetilde{X}_n(t)\Vert^p_H \nonumber\\
  \leq & \widetilde{\mathbb{E}}\liminf_{n\rightarrow\infty} \sup_{t\in[0,T]}\Vert \widetilde{X}_n(t)\Vert^p_H \nonumber\\
  \leq & \liminf_{n\rightarrow\infty}\widetilde{\mathbb{E}}\sup_{t\in[0,T]}\Vert \widetilde{X}_n(t)\Vert^p_H < \infty.
\end{align}
Similarly,
\begin{align}\label{1230.1107}
\widetilde{\mathbb{E}} \Big( \int_0^T \Vert \widetilde{X}(s)\Vert_V^{\alpha} ds \Big)^{\frac{p}{2}} < \infty.
\end{align}
%
Furthermore, by  (\ref{1230.1048}), (H4) and (H5), the following estimates hold.
\begin{lemma}\label{220120.2057}
We have the following estimates,
\begin{align}
\label{220119.1429}	& \sup_{n\in\mathbb{N}} \mathbb{E} \int_0^T \Vert A(t, \widetilde{X}_n(t)) \Vert_{V^*}^{\frac{\alpha}{\alpha-1}} dt < \infty, \\
	 & \sup_{n\in\mathbb{N}} \mathbb{E} \int_0^T \Vert P_n B(t, \widetilde{X}_n(t)) Q_n  \Vert_{L_2}^{2} dt < \infty .
\end{align}
	
\end{lemma}

The above estimates together imply that there exist $\widehat{X}\in L^{\alpha}(\widetilde{\Omega}\times[0,T], V)$, $\widetilde{\mathcal{A}}\in L^{\frac{\alpha}{\alpha-1}}(\widetilde{\Omega}\times[0,T],V^*)$ and $\widetilde{\mathcal{B}}\in L^2(\widetilde{\Omega}\times[0,T], L_2(U,H))$ such that, taking a subsequence if necessary,
\begin{gather}
\label{220123.1933}	\widetilde{X}_n \rightharpoonup \widehat{X} \quad \text{in } L^{\alpha}(\widetilde{\Omega}\times[0,T], V), \\
	\label{1230.2043} A(\cdot,\widetilde{X}_n(\cdot))\rightharpoonup \widetilde{\mathcal{A}} \quad \text{in } L^{\frac{\alpha}{\alpha-1}}(\widetilde{\Omega}\times[0,T],V^*), \\
	\label{1230.1625} P_n B(\cdot, \widetilde{X}_n(\cdot))Q_n \rightharpoonup \widetilde{\mathcal{B}} \quad \text{in } L^2(\widetilde{\Omega}\times[0,T], L_2(U,H)), \\
	\int_0^{\cdot} P_n B(s, \widetilde{X}_n(s))Q_n d \widetilde{W}(s) \rightharpoonup \int_0^{\cdot} \widetilde{\mathcal{B}}(s) d\widetilde{W}(s) \quad \text{in } L^{\infty}([0,T], L^2(\widetilde{\Omega},H)) ,
\end{gather}
here the notation ``$\rightharpoonup$" stands for the weak convergence.
\noindent Set
\begin{align} \label{1230.1633}
	\overline{X}(t):= x +\int_0^t \widetilde{\mathcal{A}}(s)ds +\int_0^t \widetilde{\mathcal{B}}(s)d\widetilde{W}(s) .
\end{align}
Then it is easy to see that
\begin{align}\label{1230.1634}
	\widetilde{X} = \widehat{X} = \overline{X}  , \quad \widetilde{\mathbb{P}}\otimes dt \text{-a.s.}.
\end{align}
In fact, the equality on the far right is known in the literature, see e.g. pages 87-88 in \cite{PR}. The first equality in (\ref{1230.1634}) follows from the uniqueness of the limits. Moreover, by Theorem 4.2.5 in \cite{LR2}, we also know that $\overline{X}$ is an $H$-valued continuous process. In view of (\ref{1230.1117}), $\widetilde{X}$ is $H$-valued, and by its continuity in $V^*$, $\widetilde{X}$ is weakly continuous in $H$. Therefore, $\widetilde{X}$ and $\overline{X}$ are indistinguishable.

\vskip 0.6cm

From now on, we will work on the new filtered probability space $(\widetilde{\Omega},\widetilde{\mathcal{F}}, \{\widetilde{\mathcal{F}}_t\}_{t\geq 0}, \widetilde{\mathbb{P}})$. However, we will drop all the superscripts $\ \widetilde{}\ $ to simplify the notations, for example, we write $\widetilde{X}_n$ and $\widetilde{X}$ as $X_n$ and $X$ respectively. Thus, (\ref{220119.1617}) reads as
\begin{align}\label{220210.1438}
\Vert X_n - X \Vert_{L^{\alpha} ([0,T],H)}+ \Vert X_n -X \Vert_{C([0,T], V^*)}\rightarrow 0  ,  \quad \mathbb{P} \text{-a.s.} .
\end{align}

\begin{lemma}\label{1230.1635}
$\mathcal{B}(\cdot)=B(\cdot, X(\cdot))$, \  $\mathbb{P}\otimes dt$ almost everywhere.
\end{lemma}
\noindent {\bf Proof}.
Since $\Vert X_n -X \Vert_{L^{\alpha} ([0,T],H)}\rightarrow 0$, \ $\mathbb{P}$-a.s., in view of (\ref{1230.1048}) and (\ref{1230.1107}), we see that
\begin{align}\label{220121.1628}
	\lim_{n\rightarrow\infty}\mathbb{E}\int_0^T \Vert X_n(t)-X(t)\Vert_H^{\kappa} dt = 0, \quad \forall\ \kappa\in[1,\alpha) .
\end{align}
Therefore, we can find a subsequence still denoted by $\{X_n\}$ such that
\begin{align}\label{1230.2211}
	\lim_{n\rightarrow\infty}\Vert X_n(t,\omega) - X(t,\omega)\Vert_H = 0 , \quad \ a.e.\ (t,\omega).
\end{align}
Now, we claim that
\begin{align}\label{1230.1616}
	\lim_{n\rightarrow\infty}\mathbb{E}\int_0^{T} \Vert P_n B(t,X_n(t))Q_n - B(t, X(t)) \Vert_{L_2}^2 dt = 0 .
\end{align}
In fact,
\begin{align}\label{231027.1610}
	& \lim_{n\rightarrow\infty}\mathbb{E}\int_0^{T} \Vert P_n B(t,X_n(t))Q_n - B(t, X(t)) \Vert_{L_2}^2 dt \nonumber\\
\leq & \lim_{n\rightarrow\infty}\mathbb{E}\int_0^{T} \Vert P_n [ B(t,X_n(t)) - B(t, X(t)) ] Q_n \Vert_{L_2}^2 dt \nonumber\\
& + \lim_{n\rightarrow\infty}\mathbb{E}\int_0^{T} \Vert P_n B(t, X(t))[ Q_n - I] \Vert_{L_2}^2 dt \nonumber\\
& + \lim_{n\rightarrow\infty}\mathbb{E}\int_0^{T} \Vert [ P_n - I ] B(t, X(t)) \Vert_{L_2}^2 dt \nonumber\\
=:& \  I_1 + I_2 + I_3.
\end{align}
By the property of Hilbert-Schmidt operators (please refer to Remark B.0.6 (iii) of [29]), we have
\begin{align}\label{231027.1541}
  I_1 \leq \lim_{n\rightarrow\infty}\mathbb{E}\int_0^{T} \Vert  B(t,X_n(t)) - B(t, X(t)) \Vert_{L_2}^2 dt = 0 ,
\end{align}
where we have used (H5), (2.55) and (2.56) in the last equality. Similarly,
\begin{align*}
  I_2 \leq & \lim_{n\rightarrow\infty}\mathbb{E}\int_0^{T} \Vert B(t, X(t))[ Q_n - I] \Vert_{L_2}^2 dt .
\end{align*}
Note that $Q_n$ is the orthogonal projection onto $\mathrm{span}\{h_1,\cdots,h_n\}$ in $U$, where $\{h_i\}_{i=1}^{\infty}$ be an orthonormal basis of Hilbert space $U$.
Then by the definition of the norm of Hilbert-Schmidt operators and the Fubini theorem,
\begin{align}\label{231027.1542}
  I_2 \leq & \lim_{n\rightarrow\infty}\mathbb{E}\int_0^{T} \sum_{i=1}^{\infty} \Vert B(t, X(t)) (Q_n - I) h_i \Vert_{H}^2 dt \nonumber\\
  = & \lim_{n\rightarrow\infty}\mathbb{E}\int_0^{T} \sum_{i>n} \Vert B(t, X(t)) h_i \Vert_{H}^2 dt \nonumber\\
  = & \lim_{n\rightarrow\infty} \sum_{i>n}\mathbb{E}\int_0^{T}  \Vert B(t, X(t)) h_i \Vert_{H}^2 dt  = 0,
\end{align}
where the zero limit is due to the convergence of the following series
\begin{align*}
 \sum_{i=1}^{\infty}\mathbb{E}\int_0^{T} \Vert B(t, X(t)) h_i \Vert_{H}^2 dt
= & \mathbb{E}\int_0^{T} \sum_{i=1}^{\infty} \Vert B(t, X(t)) h_i \Vert_{H}^2 dt \nonumber\\
= & \mathbb{E}\int_0^{T} \Vert B(t, X(t)) \Vert_{L_2}^2 dt < \infty.
\end{align*}
By the dual property of the Hilbert-Schmidt operator (please refer to Remark B.0.6 (i) of [29]), and the similar proof as (\ref{231027.1542}), we have
\begin{align}\label{231027.1605}
  I_3 = \lim_{n\rightarrow\infty}\mathbb{E}\int_0^{T} \Vert  B(t, X(t))^{*} [ P_n - I ] \Vert_{L_2(H, U)}^2 dt = 0.
\end{align}
Combining (\ref{231027.1610})-(\ref{231027.1605}) together yields (\ref{1230.1616}). The uniqueness of the limit,
(\ref{1230.1625}) and (\ref{1230.1616}) imply that $\mathcal{B}(\cdot)=B(\cdot,X(\cdot))$. 
$\blacksquare$

\vskip 0.6cm

To proceed, we will use the pseudo-monotonicity of the operator $A$. In Lemma 5.2.13 of \cite{LR2} (see also Lemma 2.2 of \cite{Liu}), $A$ was shown to be pseudo-monotone under conditions (H1), (H2) and compact embedding $V\subseteq H$.
By the similar methods, we can show that $A$ is still pseudo-monotone if we replace condition (H2) by the weaker condition (H2)$^\prime$. Here we give the proof for readers' convenience.

\begin{lemma}\label{220609.2156}
Assume (H1) and (H2)$^\prime$ hold, the embedding $V\subseteq H$ is compact. Then $A(t,\cdot)$ is pseudo-monotone from $V$ to $V^*$ for a.e. $t\in[0,T]$.  	
\end{lemma}
\noindent {\bf Proof}.
Let $\mathcal{N}\subseteq [0,T]$ be a null-set such that the assumptions (H1) and (H2) hold for any $t\in [0,T] \backslash\mathcal{N}$.
We fix $t\in [0,T]\backslash\mathcal{N}$ and denote $A(t,\cdot)$ by $A(\cdot)$.

We need to show that if  $u_n$ converges weakly to $u$ in $V$ and
	\begin{align}\label{220609.2028}
		\liminf_{n\rightarrow\infty} \langle A(u_n), u_n-u \rangle \geq 0,
	\end{align}
then for any $v\in V$,
	\begin{align}\label{220609.1543}
		\limsup_{n\rightarrow\infty} \langle A(u_n), u_n -v \rangle \leq \langle A(u), u-v \rangle .
	\end{align}
Fix $v\in V$.
Set $R = \Vert v\Vert_V + \sup_{n\in\mathbb{N}}\Vert u_n \Vert_V$.
Since the embedding $V\subseteq H$ is compact, we have $\Vert u_n - u\Vert_H \rightarrow 0$ and
\begin{align}
	\lim_{n\rightarrow\infty}\langle K_{R}(t) u_n , u_n - v \rangle = \langle K_{R}(t) u, u - v \rangle ,
\end{align}
where $K_{R}(t)$ is the constant in (H2)$^\prime$ with $R$ above.
Hence to prove (\ref{220609.1543}), it suffices to prove
\begin{align}\label{220609.2129}
	\limsup_{n\rightarrow\infty} \langle A_0(u_n), u_n -v \rangle \leq \langle A_0(u), u-v \rangle ,
\end{align}
where $A_0(u) = A(u) - K_{R}(t) u $.

By (H2)$^\prime$ and the weak convergence of
$u_n$, we have
\begin{align}
	\limsup_{n\rightarrow\infty} \langle A_0(u_n), u_n - u \rangle \leq \limsup_{n\rightarrow\infty}\langle A_0(u), u_n - u \rangle \leq 0 .
\end{align}
This together with  (\ref{220609.2028}) implies that
\begin{align}\label{220609.2124}
	\lim_{n\rightarrow\infty} \langle A_0(u_n), u_n - u \rangle = 0.
\end{align}
Let $z=u+t(v-u)$ for $t\in(0,1)$, then $\Vert z\Vert_V\leq R$. (H2)$^\prime$ yields
\begin{align}
	\langle A_0(u_n) - A_0(z), u_n - z \rangle \leq 0
\end{align}
and hence
\begin{align}
	\langle A_0(u_n), u_n - u\rangle - \langle A_0(z), u_n - u\rangle + t \langle A_0(u_n) , u - v \rangle \leq  t \langle A_0(z), u - v \rangle.
\end{align}
Taking $n\rightarrow\infty$ on both sides of the above inequality, by (\ref{220609.2124}) and the weak convergence of  $u_n$, we get
\begin{align}
	 \limsup_{n\rightarrow\infty}\langle A_0(u_n) , u_n - v \rangle \leq  \langle A_0(z), u - v \rangle.
\end{align}
Thus, letting $t\rightarrow 0$ and by the hemicontinuity (H1) we obtain (\ref{220609.2129}).
$\blacksquare$

\vskip 0.6cm

The next lemma is crucial, which shows that the following operator:
\begin{align*}
	 X(\cdot) \longmapsto A(\cdot, X(\cdot))
\end{align*}
is pseudo-monotone from $L^{\alpha}(\Omega\times[0,T], V)$ to $L^{\frac{\alpha}{\alpha-1}}(\Omega\times[0,T],V^*)$.
\begin{lemma}\label{220119.1600}
	Denote the weak convergence by the notation ``$\rightharpoonup$". If
	\begin{gather}
	X_n \rightharpoonup X \quad \text{in } L^{\alpha}(\Omega\times[0,T], V), \nonumber\\
	 A(\cdot, X_n(\cdot))\rightharpoonup \mathcal{A} \quad \text{in } L^{\frac{\alpha}{\alpha-1}}(\Omega\times[0,T],V^*), \\
\label{220119.1128}	\liminf_{n\rightarrow\infty} \mathbb{E} \int_0^T \langle A(t, X_n(t)), X_n(t) \rangle dt \geq \mathbb{E} \int_0^T \langle \mathcal{A}(t), X(t) \rangle dt ,
	\end{gather}
then $\mathcal{A}(\cdot)=A(\cdot, X(\cdot))$, \ $\mathbb{P}\otimes dt$-a.e..
\end{lemma}
\noindent {\bf Proof}. The main idea used in this proof was initiated by \cite{H}. The proof here is inspired by \cite{Liu, Sh}.

By (H3), (H4) and the Young inequality, we have
\begin{align}\label{220118.2221}
	& \langle A(t, X_n(t)), X_n(t) - X(t) \rangle \nonumber\\
	\leq & -c \Vert X_n(t) \Vert_V^{\alpha} + f(t)(1+\Vert X_n(t)\Vert_H^2) + \Vert A(t, X_n(t))\Vert_{V^*}\Vert X(t)\Vert_V \nonumber\\
	\leq & -c \Vert X_n(t) \Vert_V^{\alpha} + f(t)(1+\Vert X_n(t)\Vert_H^2) \nonumber\\
	& +  \big[f(t)+C\Vert X_n(t)\Vert_V^{\alpha}\big]^{\frac{\alpha-1}{\alpha}}\big[1+\Vert X_n(t)\Vert_H^{\beta}\big]^{\frac{\alpha-1}{\alpha}}  \Vert X(t)\Vert_V \nonumber\\
	\leq & - \frac{c}{2} \Vert X_n(t) \Vert_V^{\alpha} + f(t)(2+\Vert X_n(t)\Vert_H^2) +  C\Vert X(t)\Vert_V^{\alpha} \nonumber\\
	& +  C \Vert X_n(t)\Vert_H^{\beta(\alpha-1)} \Vert X(t)\Vert_V^{\alpha} .
\end{align}
To simplify the notation, we write
\begin{align}
	g_n(t,\omega):= &  \langle A(t, X_n(t,\omega)), X_n(t,\omega) - X(t,\omega) \rangle , \nonumber\\
	F_n(t,\omega):= & f(t)(2+\Vert X_n(t,\omega)\Vert_H^2) +  C\Vert X(t,\omega)\Vert_V^{\alpha} \nonumber\\
	& + C \Vert X_n(t,\omega)\Vert_H^{\beta(\alpha-1)} \Vert X(t,\omega)\Vert_V^{\alpha} .
\end{align}
Then (\ref{220118.2221}) reads as
\begin{align}\label{220119.1045}
	g_n(t,\omega) \leq - \frac{c}{2} \Vert X_n(t,\omega) \Vert_V^{\alpha} + F_n(t,\omega) .
\end{align}
The rest of the proof is divided into four steps.

{\bf Claim 1}: for a.e. $(t,\omega)$,  we have
\begin{align}\label{220119.1042}
	\limsup_{n\rightarrow\infty} g_n(t,\omega) \leq 0 .
\end{align}
By (\ref{1230.2211}) and Lemma \ref{220609.2156}, there exists a measurable subset $\Gamma$ of $\Omega\times[0,T]$ such that $(\Omega\times[0,T]) \backslash \Gamma $ is a $\mathbb{P}\otimes dt$-null set,
\begin{align}\label{220119.1046}
	\lim_{n\rightarrow\infty}\Vert X_n(t,\omega) - X(t,\omega)\Vert_H = 0, \quad  \forall\,  (t,\omega)\in\Gamma ,
\end{align}
and $A(t,\cdot)$ is pseudo-monotone for any $(t,\omega)\in \Gamma$.
Now take any fixed $(t,\omega)\in\Gamma$ and set
\begin{align}
	\Lambda := \{ n\in\mathbb{N}: g_n(t,\omega)> 0\} .
\end{align}
If $\Lambda$ is a finite set, then obviously (\ref{220119.1042}) holds. If $\Lambda$ is an  infinite set, then by (\ref{220119.1045}) and (\ref{220119.1046}), it follows that
\begin{align}
	\sup_{n\in\Lambda} \Vert X_n(t,\omega) \Vert_V^{\alpha} < \infty.
\end{align}
Consequently, there exists a subsequence $\{ n_i \}$ from $\Lambda$ and a element $z\in V$ such that $X_{n_i}(t,\omega)$ converges weakly to $z$ in $V$. In view of (\ref{220119.1046}), we must have $z=X(t,\omega)$ and moreover,
\begin{align}
	\lim_{\substack{n\rightarrow\infty \\ n\in\Lambda}} X_n(t,\omega) = X(t,\omega) ,
\end{align}
 weakly in $V$. Thus, using the pseudo-monotonicity of $A$ yields
\begin{align}
	\limsup_{\substack{n\rightarrow\infty \\ n\in\Lambda}} g_n(t,\omega) \leq 0.
\end{align}
On the other hand, by the definition of $\Lambda$,
\begin{align}
	\limsup_{\substack{n\rightarrow\infty \\ n\notin\Lambda}} g_n(t,\omega) \leq 0.
\end{align}
Hence Claim 1 is proved.
\vskip 0.3cm
{\bf Claim 2}:
\begin{align}
	\lim_{n\rightarrow \infty} \mathbb{E} \int_0^T g_n(t) dt = 0.
\end{align}
By (\ref{220119.1046}), we know that $F_n(t,\omega)$ converges for a.e. $(t,\omega)$. On the other hand,  it follows from (\ref{1230.1048}) that $F_n$ is uniformly integrable. Hence by a generalized Fatou Lemma (see e.g. \cite{CD}, p10), (\ref{220119.1045}) and Claim 1, we get
\begin{align}\label{220119.1130}
	\limsup_{n\rightarrow\infty} \mathbb{E} \int_0^T g_n(t) dt \leq  \mathbb{E} \int_0^T \limsup_{n\rightarrow\infty} g_n(t) dt \leq 0.
\end{align}
According to the condition (\ref{220119.1128}),
\begin{align}\label{220119.1131}
	\liminf_{n\rightarrow\infty} \mathbb{E} \int_0^T g_n(t) dt \geq 0 .
\end{align}
Thus combining (\ref{220119.1130}) and (\ref{220119.1131}) together proves Claim 2.
\vskip 0.3cm
{\bf Claim 3}: there exists a subsequence $\{ n_i\}$ such that
\begin{align}
	\lim_{i\rightarrow \infty} g_{n_i}(t,\omega) = 0 ,\quad \text{for a.e. } (t,\omega).
\end{align}
Set $g_n^{+}(t,\omega):= \max\{g_n(t,\omega), 0\}$. From Claim 1 it follows that
\begin{align}
	\lim_{n\rightarrow 0} g_{n}^{+}(t,\omega) = 0 , \quad \text{for a.e. } (t,\omega).
\end{align}
Hence by (\ref{220119.1045}) and the uniform integrability of $F_n$, we have
\begin{align}
	\lim_{n\rightarrow\infty}\mathbb{E} \int_0^T g_n^{+}(t) dt = 0.
\end{align}
Using $|g_n| = 2g_n^{+} - g_n$ and Claim 2, we see that
\begin{align}
	\lim_{n\rightarrow\infty}\mathbb{E} \int_0^T |g_n(t)|  dt = 0.
\end{align}
Claim 3 follows.
\vskip 0.3cm
{\bf Claim 4}: $\mathcal{A}(\cdot)=A(\cdot, X(\cdot))$, \ $\mathbb{P}\otimes dt$-a.e..
By (\ref{220119.1045}) and Claim 3, we have
\begin{align}\label{eq-1}
	\sup_{i\in\mathbb{N}} \Vert X_{n_i}(t,\omega) \Vert_V^{\alpha} < \infty, \quad \text{for a.e. } (t,\omega).
\end{align}
Due to (\ref{220119.1046}), we deduce from (\ref{eq-1}) that  for a.e. $(t,\omega)$, as $i\rightarrow\infty$,
\begin{align}
	X_{n_i}(t,\omega) \rightharpoonup X(t,\omega) \quad \text{weakly in } V .
\end{align}
Claim 3 together with the pseudo-monotonicity of $A$ (see Remark \ref{220119.1500}) implies that  for a.e. $(t,\omega)$
\begin{align}
	A(t,X_{n_i}(t,\omega)) \rightharpoonup A(t,X(t,\omega)) \quad \text{weakly in } V^*.
\end{align}
By (\ref{1230.2043}) and the uniqueness of the limit, we can conclude that $\mathcal{A}(\cdot)=A(\cdot, X(\cdot))$ proving  Claim 4.
$\blacksquare$

\begin{theorem}\label{220124.1213}
	There exists a probabilistically weak solution to equation (\ref{1.a}), which satisfies moment estimates (\ref{1226.2116}).
\end{theorem}
\noindent {\bf Proof}. We will show that the limit $X$ of $X_n$ obtained above is a solution to equation (\ref{1.a}). To this end, by (\ref{1230.1633}), Lemma \ref{1230.1635} and Lemma \ref{220119.1600}, we only need to verify (\ref{220119.1128}). Taking into account the equations (\ref{1230.1649}), (\ref{1230.1633}) satisfied respectively by $X_n$ and $X$, applying Ito's formula and taking expectations separately we obtain
\begin{align}
\label{220119.1637}	\mathbb{E} \Vert X_n(t) \Vert_H^2 = & \Vert P_n x\Vert_H^2 + 2 \mathbb{E} \int_0^T \langle A(t,X_n(t)), X_n(t) \rangle dt \nonumber\\
	& + \mathbb{E} \int_0^T \Vert P_n B(t,X_n(t)) Q_n \Vert_{L_2}^2 dt , \\
\label{220119.1638}	\mathbb{E} \Vert X(t) \Vert_H^2 = & \Vert  x\Vert_H^2 + 2 \mathbb{E} \int_0^T \langle \mathcal{A}(t), X(t) \rangle dt \nonumber\\
	& + \mathbb{E} \int_0^T \Vert B(t,X(t)) \Vert_{L_2}^2 dt .
\end{align}
Since $\Vert X_n - X \Vert_{C([0,T],V^*)} \rightarrow 0$ (see (\ref{220210.1438})), by the lower semicontinuity of $\Vert \cdot \Vert_H$ in $V^*$ and Fatou's lemma, we have
\begin{align}
	\mathbb{E} \Vert X(t) \Vert_H^2 \leq \mathbb{E} \liminf_{n\rightarrow\infty}\Vert X_n(t) \Vert_H^2 \leq \liminf_{n\rightarrow\infty} \mathbb{E} \Vert X_n(t) \Vert_H^2.
\end{align}
Hence in view of (\ref{1230.1616}), and comparing (\ref{220119.1637}) and (\ref{220119.1638}), we see that (\ref{220119.1128}) holds.
Moreover, the moment estimates (\ref{1226.2116}) for $X$ follow from the  estimates (\ref{1230.1117}) and (\ref{1230.1107}).
$\blacksquare$

\begin{theorem}\label{1231.0755}
If (H2) is satisfied, then the pathwise uniqueness holds for solutions of equation (\ref{1.a}).
\end{theorem}
\noindent {\bf Proof}.  Let $X$ and $X^{\prime}$ be two solutions of equation (\ref{1.a}) defined on a same probability space and driven by the same Wiener process, with initial values $X(0)=x$ and $X^{\prime}(0)=x^{\prime}$ respectively. Set
\begin{align}
	\varphi(t):=\exp\left(-\int_0^t \big[f(r)+\rho(X(r))+\eta(X^{\prime}(r))\big]dr \right).
\end{align}
Then $\varphi$ is a continuous process of finite variation. By Ito's formula and (H2), we have for any $t\in [0,T]$,
\begin{align}\label{220119.1739}
	& \varphi(t)\Vert X(t)-X^{\prime}(t)\Vert_H^2 \nonumber\\
	= & \Vert x -x^{\prime}\Vert_H^2 + \int_0^t \varphi(s) \Big\{2\langle A(s,X(s))-A(s,X^{\prime}(s)), X(s)-X^{\prime}(s)\rangle \nonumber\\
	& + \Vert B(s,X(s)) -B(s,X^{\prime}(s))\Vert_{L_2}^2  \nonumber\\
	& - \big[f(s)+\rho(X(s))+\eta(X^{\prime}(s))\big]\Vert X(s)-X^{\prime}(s)\Vert_H^2  \Big\} ds \nonumber\\
	& + 2\int_0^T \varphi(s) ( X(s)-X^{\prime}(s), \big[B(s,X(s)) - B(s,X^{\prime}(s))\big]dW(s) ) \nonumber\\
	\leq &  \Vert x -x^{\prime}\Vert_H^2 + 2\int_0^t \varphi(s) ( X(s)-X^{\prime}(s), \big[B(s,X(s)) - B(s,X^{\prime}(s))\big]dW(s) ) .
\end{align}
Let $\{\sigma_l\}\uparrow\infty$ be a sequence of stopping times such that the local martingale in the above inequality is a martingale. Then taking the expectation on both sides of the above inequality, we get
\begin{align}\label{1231.1341}
	 \mathbb{E}\Big[ \varphi(t\wedge\sigma_l)\Vert X(t\wedge\sigma_l)-X^{\prime}(t\wedge\sigma_l)\Vert_H^2 \Big] \leq \Vert x - x^{\prime}\Vert_H^2 .
\end{align}
Letting $l\rightarrow\infty$ and applying Fatou's lemma yield that for any $t\in [0,T]$,
\begin{align}\label{1231.1342}
	\mathbb{E} \Big[\varphi(t)\Vert X(t)-X^{\prime}(t)\Vert_H^2 \Big] \leq \Vert x -x^{\prime}\Vert_H^2 .
\end{align}
By (H2) and the estimate (\ref{1226.2116}), we see that
\begin{align}
	\int_0^T \big[f(r)+\rho(X(r))+\eta(X^{\prime}(r))\big]dr < \infty, \quad \mathbb{P}\text{-a.s. in } \Omega.
\end{align}
According to the definition of $\varphi(t)$ and the above estimate, for any $t\in [0,T]$ we have $\varphi(t,\omega) > 0$ for $\mathbb{P}$-a.s. $\omega\in\Omega$. Hence (\ref{1231.1342}) and the continuity of of $X$ and $X^{\prime}$ in $H$  in particular imply the pathwise uniqueness of solutions to equation (\ref{1.a}).
$\blacksquare$

\vskip 0.6cm
Theorem \ref{1227.2229} is a combination of the above Theorem \ref{220124.1213} and Theorem \ref{1231.0755}. Next we give
\vskip 0.4cm
\noindent {\bf Proof of Theorem \ref{220119.1651}}. For $M>0$, we define the stopping time
\begin{align*}
	\sigma_{n}^M := & \inf\Big\{t\in [0,T]: \Vert X(t, x_n)\Vert_H> M \Big\} \nonumber\\
	 & \wedge\inf\Big\{t\in [0,T]: \int_0^t\Vert X(s, x_n)\Vert_{V}^{\alpha}ds >M\Big\} \nonumber\\
	 & \wedge \inf\Big\{t\in [0,T]: \Vert X(t, x)\Vert_H> M \Big\} \nonumber\\
	 & \wedge\inf\Big\{t\in [0,T]: \int_0^t\Vert X(s, x)\Vert_{V}^{\alpha}ds >M\Big\} \wedge T
\end{align*}
with the convention $\inf\emptyset = +\infty$. Then by the moment estimates (\ref{1226.2116}) for  the solutions we have
\begin{align}\label{1231.1400}
	\lim_{M\rightarrow\infty}\sup_{n\in \mathbb{N}} \mathbb{P}(\sigma_{n}^M <T) =0 .
\end{align}
From (\ref{1231.1341}) and (\ref{1231.1342}), it follows that
\begin{align}
	 \mathbb{E}\Big[\varphi_n(t\wedge\sigma_n^M)\Vert X(t\wedge\sigma_n^M, x_n) -X(t\wedge\sigma_n^M, x)\Vert_H^2 \Big]
	\leq  \Vert x_n - x \Vert_H^2 ,
\end{align}
where
\begin{align}
	\varphi_n(t):=\exp\left(-\int_0^t \big[f(r)+\rho(X(r,x_n))+\eta(X(r,x))\big]dr \right).
\end{align}
Now, for any $\epsilon>0$, there exists a constant $C_M>0$ such that
\begin{align}\label{220119.2002}
	& \mathbb{P}\big(\Vert X(t,x_n)-X(t,x)\Vert_H \geq\epsilon\big) \nonumber\\
	\leq & \mathbb{P}\big(\Vert X(t,x_n)-X(t,x)\Vert_H \geq\epsilon, \sigma^n_M\geq T\big)
	 + \mathbb{P}(\sigma^n_M <T) \nonumber\\
	\leq & \frac{1}{\epsilon^2 C_M} \mathbb{E} \big[ \varphi_n(t\wedge\sigma_n^M)\Vert X(t\wedge\sigma_n^M, x_n)-X(t\wedge\sigma_n^M, x)\Vert_H^{2} \big ]  + \mathbb{P}(\sigma^n_M <T) \nonumber\\
	\leq & \frac{1}{\epsilon^2 C_M} \Vert x_n - x\Vert_H^2 + \sup_{n\in\mathbb{N}}\mathbb{P}(\tau^n_M <T) .
\end{align}
In view of (\ref{1231.1400}), we let $n\rightarrow\infty$ and then  $M\rightarrow\infty$ to get that for any $t\in[0,T]$,
\begin{align}
	\lim_{n\rightarrow\infty}\Vert X(t, x_n)-X(t, x)\Vert_H = 0,  \quad \text{in probability } \ \mathbb{P} .
\end{align}
Furthermore, by Lemma \ref{1231.1011}, we see that for any $p\geq 2$,
	\begin{align}\label{220119.1942}
		\sup_{n\in\mathbb{N}} \mathbb{E} \Big[ \sup_{t\in [0,T]}\Vert X(t, x_n)\Vert_H^p \Big]   <\infty .
	\end{align}	
Hence it follows that
\begin{align}
	\lim_{n\rightarrow\infty}\mathbb{E}\int_0^T \Vert X(t,x_n) - X(t,x) \Vert_H^2 dt = 0 .
\end{align}
In particular,
\begin{align}\label{220119.1941}
	\Vert X(t,x_n) - X(t,x) \Vert_H \xrightarrow[n\rightarrow\infty]{} 0 \quad \text{in measure }\  \mathbb{P}\otimes dt .
\end{align}
Therefore, by (H5) and (\ref{220119.1942}) we get
\begin{align}\label{220119.1947}
	\lim_{n\rightarrow\infty}\mathbb{E} \int_0^T \Vert B(t,X(t,x_n)) - B(t,X(t,x)) \Vert_{L_2}^2 dt = 0 .
\end{align}
Now, by (\ref{220119.1739}), the BDG inequality and the Young inequality, we have
\begin{align}\label{eq-2}
	& \mathbb{E} \sup_{t\leq T\wedge\sigma_n^M} \Big[ \varphi_n(t) \Vert X(t,x_n) - X(t,x) \Vert_H^2 \Big] \nonumber\\
	\leq & \Vert x_n -x \Vert_H^2 + 2\mathbb{E} \sup_{t\leq T\wedge\sigma_n^M}\Big| \int_0^t \varphi_n(s) ( X(s,x_n)-X(s,x), \nonumber\\
	&~~~~~~~~~~~~~~~~~~~~~\big[B(s,X(s,x_n)) - B(s,X(s,x))\big]dW(s) ) \Big| \nonumber\\
	\leq & \Vert x_n -x \Vert_H^2 + C\mathbb{E} \Big( \int_0^{T\wedge\sigma_n^M} \varphi_n(t)^2 \Vert X(t,x_n) - X(t,x) \Vert_H^2 \nonumber\\
	& ~~~~~~~~~~~~~~~~~~\times \Vert B(t,X(t,x_n)) - B(t,X(t,x)) \Vert_{L_2}^2 dt \Big)^{\frac{1}{2}}  \nonumber\\
	\leq & \Vert x_n -x \Vert_H^2 + \frac{1}{2} \mathbb{E} \sup_{t\leq T\wedge\sigma_n^M} \Big[ \varphi_n(t) \Vert X(t,x_n) - X(t,x) \Vert_H^2 \Big] \nonumber\\
	&  + C \mathbb{E} \int_0^{T\wedge\sigma_n^M} \varphi_n(t) \Vert B(t,X(t,x_n)) - B(t,X(t,x)) \Vert_{L_2}^2 dt .
\end{align}
 (\ref{eq-2}) and (\ref{220119.1947}) imply
\begin{align}
	\lim_{n\rightarrow\infty}\mathbb{E} \sup_{t\leq T\wedge\sigma_n^M} \Big[ \varphi_n(t) \Vert X(t,x_n) - X(t,x) \Vert_H^2 \Big] = 0.
\end{align}
Arguing as (\ref{220119.2002}) again yields
\begin{align}\label{220121.2337}
	\lim_{n\rightarrow\infty}\sup_{t\in [0,T]}\Vert X(t, x_n)-X(t, x)\Vert_H = 0,  \quad \text{in probability } \ \mathbb{P} .
\end{align}
Hence it follows from (\ref{220119.1942}) that
\begin{align}
	\lim_{n\rightarrow\infty}\mathbb{E}\Big[ \sup_{t\in [0,T]}\Vert X(t, x_n)-X(t, x)\Vert_H^p \Big] = 0,
\end{align}
completing the proof. $\blacksquare$

\section{Part II}\label{220123.1616}
\setcounter{equation}{0}

In this part, we will allow the dependence of  $\Vert B(t,u) \Vert_{L_2}$ on the $V$-norm $\Vert u \Vert_{V}$. In the situation of classical stochastic partial differential equations, this typically means that  $B(t,u)$ is allowed to depend also on the gradient $\nabla u$ of the solution function $u$. We will modify the arguments used in Section \ref{220406.2023} to establish the well-posedness of equation (\ref{1.a}) under a new set of local monotone conditions which are slight adjustment of the hypotheses in Section \ref{220406.2023}. Let us now introduce the assumptions.

\vskip 0.6cm

Let $f\in L^1([0,T],\mathbb{R}_{+})$, $\alpha\in (1,\infty)$ and $\beta \in [0,\infty)$.
\begin{itemize}
 	\item [(H2)*] There exist nonnegative constants $\theta\in [0,\alpha), \gamma, \lambda $ and $C$ such that for a.e. $t\in [0,T]$, and any $u,v\in V$,
\begin{align}\label{220124.1150}
& 2 \langle A(t,u)-A(t,v), u-v \rangle + \Vert B(t,u)-B(t,v)\Vert_{L_2}^2 \nonumber\\
\leq & [f(t)+ \rho(u)+\eta(v) ] \Vert u-v\Vert_{H}^2 ,
\end{align}
where $\rho$ and $\eta$ are two measurable functions from $V$ to $\mathbb{R}$ satisfying
\begin{align}
\label{220120.1632}	|\rho(u)| \leq &  C (1+\Vert u\Vert_H^{\lambda}) + C \Vert u\Vert_V^{\theta} (1+\Vert u\Vert_H^{\gamma}), \\
\label{220121.1756}	|\eta(u)| \leq & C(1+\Vert u \Vert_H^{2+\beta})+ C\Vert u\Vert_V^{\alpha}(1+\Vert u\Vert_H^{\beta}).
\end{align}


%

	\item [(H3)*] There exists a constant $L_A>0$ such that for a.e. $t\in [0,T]$, and any $u\in V$,
		\begin{align}
			\langle A(t,u),u \rangle \leq f(t) (1+\Vert u\Vert_{H}^2) - L_A \Vert u \Vert_{V}^{\alpha} .
		\end{align}
	\item [(H4)*] There exists nonnegative constant $C$ such that for a.e. $t\in [0,T]$, and any $u\in V$,
		\begin{align}
			\Vert A(t,u)\Vert_{V^*}^{\frac{\alpha}{\alpha-1}} \leq f(t) ( 1+ \Vert u \Vert_H^{2+\beta}) + C \Vert u\Vert_V^{\alpha}(1+\Vert u\Vert_H^{\beta}) .
		\end{align}
	\item [(H5)*] There exists $g \in L^{1}([0,T],\mathbb{R}_+)$ and a constant $L_B\geq 0$
such that for a.e. $t\in [0,T]$, and any $u\in V$,
\begin{align}\label{220602.2309}
	\Vert B(t,u)\Vert^2_{L_2} \leq g(t) (1+ \Vert u\Vert_{H}^2 ) + L_B \Vert u \Vert_V^{\alpha}.
\end{align}

\end{itemize}
\begin{remark}
 The stronger condition $\theta<\alpha$ (than that in (H2)) in (H2)* is important to the proof of Theorem \ref{220120.1643} below. As the positions of $\rho$ and $\eta$ in (\ref{220124.1150}) are symmetric,  $\rho$ and $\eta$ can interchange in (\ref{220120.1632}) and (\ref{220121.1756}).  In contrast to (H5) of Section 2, in (H5)* there is no assumption of continuity of $B$ with respect to $H$-norm and $B$ can depend on the $V$-norm, which is the main focus of this section.
\end{remark}

The main result in this part reads as follows

\begin{theorem}\label{220120.1643}
	Suppose that the embedding $V\subseteq H$ is compact and that (H1), (H2)*, (H3)*, (H4)*, (H5)* hold with
	\begin{align}\label{231028.2048}
		L_B < \frac{2L_A}{\chi} ,
	\end{align}
where
\begin{align}\label{220122.2146}
\chi =
\begin{cases}
	\max\{1+\beta, 1+\lambda, 1+ \gamma+ \frac{2\theta}{\alpha} \}, \quad & \text{when } \alpha \leq 2 , \\
	\max\{1+\beta, 3+\lambda-\alpha, 3+\gamma+\theta-\alpha\}, \quad & \text{when } \alpha > 2 .
\end{cases}
\end{align}
Then for any initial value $x\in H$, there exists a unique probabilistically strong solution to equation (\ref{1.a}). Furthermore, for any
\begin{align}\label{220120.2002}
	2\leq p < 1+ \frac{2L_A}{L_B} ,
\end{align}
we have the following moment estimate,
\begin{align}\label{220121.2344}
	\mathbb{E} \Big\{\sup_{t\in [0,T]} \Vert X(t)\Vert_H^p\Big\} +  \mathbb{E}\left\{ \Big( \int_0^T \Vert X(t)\Vert_V^{\alpha} dt \Big)^{\frac{p}{2}}  \right\} < \infty.
\end{align}
Moreover, let $\{ x_n \}_{n=1}^{\infty}$ and $x$ be a sequence in $H$ with $ \Vert x_n -x \Vert_H \rightarrow 0$, and let $X(t,x)$ be the unique solution of (\ref{1.a}) with the initial value $x$.
Then
\begin{align}\label{220123.2050}
	 \lim_{n\rightarrow\infty}\mathbb{E}\Big[ \sup_{t\in [0,T]}\Vert X(t, x_n) - X(t, x) \Vert_H^p \Big] = 0 ,
\end{align}
for $p$ satisfying (\ref{220120.2002}).
\end{theorem}

\vskip 0.6cm

The rest of this part is devoted to the proof of Theorem \ref{220120.1643}.
We will assume the conditions of Theorem \ref{220120.1643} throughout.

As in Section \ref{220406.2023}, our starting point is the sequence of Galerkin approximating solutions.
The global existence of probabilistically weak solutions to the finite dimensional Galerkin approximating stochastic differential equations (\ref{1228.2225}) follows by the same reason as in Part I.
Since we do not assume that $B$ is  continuous on $H$, some of the proofs (e.g. the proof of Theorem \ref{220124.1213}) are not valid. The pseudo-monotonicity argument doesn't work in this case. We will instead combine the tightness of the Galerkin approximations with the monotonicity argument.

\vskip 0.6cm

Now we establish the uniform moment estimates of Galerkin approximating solutions $\{Y_n\}$ under the new assumptions. Since we will pass to a new probability space as in Section \ref{220406.2023}, for the simplicity of notations, we write $Y_n$ as $X_n$.
\begin{lemma}\label{220121.1807}
	For any $p$ satisfying (\ref{220120.2002}), there exists a constant $C_p$ such that
	\begin{align}\label{220120.2008}
		& \sup_{n\in\mathbb{N}} \bigg\{ \mathbb{E} \Big[ \sup_{t\in [0,T]} \Vert X_n(t)\Vert_H^p \Big] + \mathbb{E}\int_0^T \Vert X_n(t) \Vert_V^{\alpha} \Vert X_n(t)\Vert_H^{p-2} dt \nonumber\\
		+ & \mathbb{E} \Big( \int_0^T \Vert X_n(t)\Vert_V^{\alpha} dt \Big)^{\frac{p}{2}}\bigg\} \leq C_p(1+\Vert x \Vert_H^p).
	\end{align}	

\end{lemma}
\noindent {\bf Proof}. Using Ito's formula it follows that
\begin{align}\label{220120.1710}
	\Vert X_n(t) \Vert_H^p \leq & \Vert P_n x \Vert_H^p + p\int_0^t \Vert X_n(s) \Vert_H^{p-2}  \langle A(s, X_n(s)), X_n(s) \rangle  ds \nonumber\\
	& + \frac{p(p-1)}{2} \int_0^t \Vert X_n(s) \Vert_H^{p-2} \Vert B(s,X_n(s)) \Vert_{L_2}^2 ds \nonumber\\
	& + p\int_0^t \Vert X_n(s) \Vert_H^{p-2} ( X_n(s), B(s,X_n(s))Q_n dW(s) ) \nonumber\\
	\leq & \Vert P_n x \Vert_H^p + p\int_0^t \Vert X_n(s) \Vert_H^{p-2}  \big[ -L_A \Vert X_n(s) \Vert_V^{\alpha} + f(s)(1+\Vert X_n(s) \Vert_H^2 ) \big]  ds \nonumber\\
	& + \frac{p(p-1)}{2} \int_0^t \Vert X_n(s) \Vert_H^{p-2} \big[ L_B \Vert X_n(s)\Vert_V^{\alpha} + g(s)(1+\Vert X_n(s)\Vert_H^2) \big] ds \nonumber\\
	& + p\int_0^t \Vert X_n(s) \Vert_H^{p-2} ( X_n(s), B(s,X_n(s))Q_n dW(s) )  .
\end{align}
Rearranging terms and using stopping arguments, we get
\begin{align}\label{eq-3}
	& \mathbb{E}\Vert X_n(t) \Vert_H^p + p\Big(L_A - \frac{p-1}{2} L_B \Big) \mathbb{E}\int_0^t \Vert X_n(s) \Vert_V^{\alpha} \Vert X_n(s)\Vert_H^{p-2} ds \nonumber\\
	\leq & \Vert x\Vert_H^p + C_p \int_0^T \big[f(s)+g(s)\big] ds + C_p \mathbb{E}\int_0^t \big[f(s)+g(s)\big] \Vert X_n(s) \Vert_H^p ds .
\end{align}
The range of the parameter $p$ implies
\begin{align}
		 L_A -\frac{p-1}{2}L_B > 0 .
\end{align}
By (\ref{eq-3}) and  Gronwall's inequality we obtain
\begin{align}\label{220120.1735}
	\sup_{n\in\mathbb{N}} \left\{ \sup_{t\in [0,T]}\mathbb{E}\Vert X_n(t) \Vert_H^p + \mathbb{E}\int_0^T \Vert X_n(s) \Vert_V^{\alpha} \Vert X_n(s)\Vert_H^{p-2} ds  \right\}\leq  C_p(1+\Vert x \Vert_H^p) .
\end{align}
Again using (\ref{220120.1710}) and (H3)*, we have
\begin{align}\label{220120.1727}
	& \mathbb{E}\Big[ \sup_{t\in [0,T]} \Vert X_n(t) \Vert_H^p \Big] + p L_A \mathbb{E}\int_0^T \Vert X_n(s) \Vert_V^{\alpha} \Vert X_n(s) \Vert_H^{p-2} ds \nonumber\\
	\leq & \Vert P_n x \Vert_H^p + p \mathbb{E}\int_0^T \Vert X_n(s) \Vert_H^{p-2} f(s)(1+\Vert X_n(s)\Vert_H^2 )  ds \nonumber\\
	& + \frac{p(p-1)}{2} \mathbb{E}\int_0^T \Vert X_n(s) \Vert_H^{p-2} \Vert B(s,X_n(s)) \Vert_{L_2}^2 ds \nonumber\\
	& + p \mathbb{E} \sup_{t\in [0,T]}\left| \int_0^t \Vert X_n(s) \Vert_H^{p-2} ( X_n(s), B(s,X_n(s))Q_n dW(s) ) \right| .
\end{align}
Similarly to (\ref{220117.1435}), by the BDG and Young inequalities we deduce that
\begin{align}\label{220120.1728}
	& p \mathbb{E} \sup_{t\in [0,T]}\left| \int_0^t \Vert X_n(s) \Vert_H^{p-2} ( X_n(s), B(s,X_n(s))Q_n dW(s) ) \right| \nonumber\\
	\leq & \frac{1}{2} \mathbb{E}\sup_{t\in [0,T]} \Vert X_n(t) \Vert_H^p  + C_p \mathbb{E}\int_0^T \Vert X_n(s) \Vert_H^{p-2} \Vert B(s,X_n(s)) \Vert_{L_2}^2 ds .
\end{align}
Combining (\ref{220120.1727}) and (\ref{220120.1728}) together and using (H5)*, we otain
\begin{align}
	& \mathbb{E}\Big[ \sup_{t\in [0,T]} \Vert X_n(t) \Vert_H^p \Big] \nonumber\\
	\leq & 2\Vert x \Vert_H^p + C_p\mathbb{E}\int_0^T \big[ f(s)+g(s)\big]ds \times \Big( 1+ \sup_{s\in [0,T]} \mathbb{E} \Vert X_n(s) \Vert_H^p \Big) \nonumber\\
	&  + C_p\mathbb{E}\int_0^T \Vert X_n(s)\Vert_V^{\alpha}\Vert X_n(s)\Vert_H^{p-2} ds .
\end{align}
Therefore, it follows from (\ref{220120.1735}) that
\begin{align}\label{220120.1812}
	\sup_{n\in\mathbb{N}} \mathbb{E} \Big[ \sup_{t\in [0,T]} \Vert X_n(t) \Vert_H^p \Big] \leq C_p (1+\Vert x\Vert_H^p) .
\end{align}
By (\ref{220117.1455}), (H3)* and (H5)*, we have
\begin{align}\label{220120.1810}
	& \mathbb{E}\Big(\int_0^T \Vert X_n(t)\Vert_V^{\alpha} dt \Big)^{\frac{p}{2}} \nonumber\\
	\leq & C_p \Vert P_n x \Vert_H^p + C_p \mathbb{E}\left| \int_0^T \big[ f(t)+g(t)\big]\big( 1+ \Vert X_n(t)\Vert_H^2 \big) dt \right|^{\frac{p}{2}} \nonumber\\
	& + C_p \mathbb{E}\left| \int_0^T ( X_n(t), B(t,X_n(t))Q_n dW(t) ) \right|^{\frac{p}{2}}  .
\end{align}
By the BDG inequality, (H5)* and the Young inequality, we get
\begin{align}\label{220120.1811}
	& C_p \mathbb{E}\left| \int_0^T ( X_n(t), B(t,X_n(t))Q_n dW(t) ) \right|^{\frac{p}{2}} \nonumber\\
	\leq & C_p \mathbb{E} \Big( \int_0^T \Vert X_n(t)\Vert_H^2 \Big[ g(t) (1+\Vert X_n(t)\Vert_H^2 ) + L_B \Vert X_n(t)\Vert_V^{\alpha} \Big] dt \Big)^{\frac{p}{4}} \nonumber\\
	\leq & C_p \mathbb{E} \bigg\{ \Big[ 1+ \sup_{t\in [0,T]} \Vert X_n(t) \Vert_H^4 \Big] \int_0^T g(t) dt \bigg\}^{\frac{p}{4}} \nonumber\\
	& + C_p \mathbb{E} \bigg\{ \sup_{t\in [0,T]} \Vert X_n(t) \Vert_H^2 \cdot \int_0^T \Vert X_n(t)\Vert_V^{\alpha} dt \bigg\}^{\frac{p}{4}} \nonumber\\
	\leq & \frac{1}{2} \mathbb{E} \Big( \int_0^T \Vert X_n(t) \Vert_V^{\alpha} dt  \Big)^{\frac{p}{2}} + C_p \Big(1+ \mathbb{E} \Big[ \sup_{t\in [0,T]} \Vert X_n(t) \Vert_H^p \Big] \Big) .
\end{align}
Combining (\ref{220120.1810}) and (\ref{220120.1811}) together, we derive that
\begin{align}\label{eq-4}
	\mathbb{E}\Big(\int_0^T \Vert X_n(t)\Vert_V^{\alpha} dt \Big)^{\frac{p}{2}} \leq C_p \Big(1+ \Vert x \Vert_H^p + \mathbb{E} \Big[ \sup_{t\in [0,T]} \Vert X_n(t) \Vert_H^p \Big] \Big) .
\end{align}
(\ref{eq-4}) and (\ref{220120.1812}) together gives the desired estimate (\ref{220120.2008}). We complete the proof of this lemma.
$\blacksquare$

\vskip 0.6cm

Repeating the proof of Lemma \ref{1229.2041}, we see that the family of the laws of $\{X_n\}_{n=1}^{\infty}$ is tight in the space $L^{\alpha}([0,T],H)$. We would like to point out that the laws of $\{X_n\}_{n=1}^{\infty}$  might not be tight in the space $C([0,T], V^*)$, since according to the proof of Lemma \ref{1229.2041}, (\ref{1229.1651}) may not hold under (H5)*. Thus by the Prohorov theorem and the modified Skorokhod representation theorem, we can pass to a new filtered probability space (still written as $(\Omega,\mathcal{F},\{\mathcal{F}_t\},\mathbb{P})$) similarly as in Section 2, and we may as well  assume there exists an $\{\mathcal{F}_t\}$-adapted process $X$ such that
\begin{align}
	\Vert X_n - X \Vert_{L^{\alpha}([0,T],H)} \rightarrow 0, \quad \mathbb{P}\text{-a.s.}.
\end{align}
By (H4)*, (H5)* and Lemma \ref{220121.1807},
the following uniform estimates hold,
\begin{align}
\label{220122.2047}  & \sup_{n\in\mathbb{N}} \mathbb{E} \int_0^T \Vert A(t, X_n(t)) \Vert_{V^*}^{\frac{\alpha}{\alpha-1}} dt < \infty, \\
\label{220122.2048}	 & \sup_{n\in\mathbb{N}} \mathbb{E} \int_0^T \Vert P_n B(t, X_n(t)) Q_n \Vert_{L_2}^{2} dt < \infty .
\end{align}
Similarly to (\ref{220121.1628}) and (\ref{1230.2211}), we have (take a subsequence if necessary)
\begin{align}\label{220121.1753}
	\lim_{n\rightarrow\infty}\Vert X_n(t,\omega) - X(t,\omega)\Vert_H = 0 , \quad \text{a.e.}\ (t,\omega).
\end{align}
Also as in Section \ref{220406.2023}, $X_n$ converges weakly to $X$ in $L^{\alpha}(\Omega\times[0,T], V)$, and
\begin{align}\label{220121.2329}
	X(t) = x + \int_0^t \mathcal{A}(s) ds + \int_0^t \mathcal{B}(s) dW(s) , \quad \mathbb{P}\otimes dt \text{-a.e.},
\end{align}
where $\mathcal{A}$ and $\mathcal{B}$ are limits of the following weak convergence (up to a subsequence),
\begin{align}
\label{eq-6} A(\cdot, X_n(\cdot)) &\rightharpoonup \mathcal{A} \quad \text{in } L^{\frac{\alpha}{\alpha-1}}(\Omega\times[0,T],V^*), \\
 \label{eq-7} P_n B(\cdot, X_n(\cdot))Q_n &\rightharpoonup \mathcal{B} \quad \text{in } L^2(\Omega\times[0,T], L_2(U,H)).
\end{align}
The next result concludes that $X$ is a solution to the equation (\ref{1.a}).

\begin{lemma}\label{220121.2330}
$\mathcal{B}(\cdot) = B(\cdot, X(\cdot))$ and  $\mathcal{A}(\cdot) = A(\cdot, X(\cdot))$, $\mathbb{P}\otimes dt$-a.e..
Consequently, $X$ is a solution to equation (\ref{1.a}).
 \end{lemma}
\noindent {\bf Proof}. Fix any $T>0$,  $dt$ denotes the Lebesgue measure on the interval $[0, T]$.
Let $u$ be any given $H$-valued continuous adapted process such that
\begin{align}\label{220122.1716}
	\mathbb{E}\Big[\sup_{t\in [0,T]}\Vert u(t)\Vert_H^{2+\beta}\Big] + \mathbb{E}\int_0^T \Vert u(t) \Vert_V^{\alpha} \big( 1+ \Vert u(t) \Vert_H^{\beta}\big) dt < \infty .
\end{align}
For $M>0$, we define the stopping time
\begin{align}
	\tau_{u}^M := &  \inf \{t\in [0,T]: \Vert u(t)\Vert_H^2 > M \} \nonumber\\
  & \wedge \inf\Big\{t\in [0,T]: \int_0^t\Vert u(s)\Vert_{V}^{\alpha}ds > M\Big\} \wedge T
\end{align}
with the convention $\inf\emptyset = +\infty$. Then
\begin{align}\label{220122.1736}
	\lim_{M\rightarrow\infty} \mathbb{P}(\tau_u^M < T) = 0.
\end{align}
For any $\epsilon > 0$,
\begin{align}
	& \mathbb{P}\otimes dt \Big( \big\{ (t,\omega): \Vert X_n(t\wedge\tau_u^M(\omega),\omega) - X(t\wedge \tau_u^M(\omega),\omega) \Vert_H > \epsilon \big\} \Big) \nonumber\\
	\leq & \mathbb{P}\otimes dt \Big( \big\{ (t,\omega): \Vert X_n(t\wedge\tau_u^M(\omega),\omega) - X(t\wedge \tau_u^M(\omega),\omega) \Vert_H > \epsilon \big\} \nonumber\\
	& \cap \big\{(t,\omega): \tau_u^M \geq T \big\}\Big)  +  \mathbb{P}\otimes dt \Big( \big\{(t,\omega): \tau_u^M(\omega) < T \big\}\Big) \nonumber\\
	\leq & \mathbb{P}\otimes dt \Big( \big\{ (t,\omega): \Vert X_n(t,\omega) - X(t,\omega) \Vert_H > \epsilon \big\} \Big)  \nonumber\\
	& + T \mathbb{P} \Big( \big\{ \omega : \tau_u^M(\omega) < T \big\}\Big)  .
\end{align}
Letting $n\rightarrow\infty$ and $M\rightarrow\infty$, in view of (\ref{220121.1753}) and (\ref{220122.1736}), we obtain
\begin{align}\label{eq-5}
	\lim_{\substack{n\rightarrow\infty \\ M\rightarrow\infty}}\Vert X_n(t\wedge\tau_u^M) - X(t\wedge \tau_u^M ) \Vert_H = 0 , \quad \text{in measure } \mathbb{P}\otimes dt .
\end{align}
Hence for any $\psi\in L^{\infty}([0,T],\mathbb{R}_{+})$, by Lemma \ref{220121.1807} and (\ref{eq-5}) it follows that
\begin{align}\label{220122.2055}
	& \lim_{M\rightarrow\infty}\mathbb{E} \int_0^T \psi(t) \Big[ \Vert X (t\wedge\tau_u^M) \Vert_H^2 - \Vert  x \Vert_H^2 \Big]  dt  \nonumber\\
	= & \lim_{M\rightarrow\infty}\liminf_{n\rightarrow\infty} \mathbb{E} \int_0^T \psi(t) \Big[ \Vert X_n (t\wedge\tau_u^M) \Vert_H^2 - \Vert P_n x \Vert_H^2 \Big]  dt .
\end{align}
Using Ito's formula and inserting terms we get
\begin{align}
	& \mathbb{E} \Vert X_n (t\wedge\tau_u^M) \Vert_H^2 - \Vert P_n x \Vert_H^2  \nonumber\\
	= & \mathbb{E} \int_0^{t\wedge\tau_{u}^M} \Big[ 2\langle A(s, X_n(s)), X_n(s) \rangle + \Vert P_n B(s,X_n(s)) Q_n \Vert_{L_2}^2 \Big] ds \nonumber\\
	\leq & \mathbb{E} \int_0^{t\wedge\tau_{u}^M} \Big[ 2\langle A(s, X_n(s)) - A(s,u(s)), X_n(s) - u(s)\rangle \nonumber\\
	& + \Vert B(s,X_n(s)) -B(s,u(s)) \Vert_{L_2}^2  \Big] ds \nonumber\\
	& +  \mathbb{E} \int_0^{t\wedge\tau_{u}^M} \Big[ 2\langle A(s, X_n(s)), u(s)\rangle + 2\langle A(s,u(s)), X_n(s) \rangle - 2\langle A(s,u(s)), u(s) \rangle \nonumber\\
	& + 2 \big( B(s,X_n(s)), B(s,u(s)) \big)_{L_2} - \Vert B(s,u(s)) \Vert_{L_2}^2 \Big] ds .
\end{align}
Hence by the Fubini theorem, (H2)*, Lemma \ref{220121.1807}, (\ref{eq-6}) and (\ref{eq-7}), it follows that
\begin{align}\label{220122.2056}
	& \liminf_{n\rightarrow\infty} \mathbb{E} \int_0^T \psi(t) \Big[ \Vert X_n (t\wedge\tau_u^M) \Vert_H^2 - \Vert P_n x \Vert_H^2 \Big]  dt \nonumber\\
	\leq & \liminf_{n\rightarrow\infty} \mathbb{E} \int_0^T \psi(t) \int_0^{t\wedge\tau_{u}^M} \big[ f(s) + \rho(X_n(s)) + \eta(u(s)) \big] \Vert X_n(s) - u(s) \Vert_H^2 ds  dt \nonumber\\
	& + \mathbb{E} \int_0^T \psi(t) \int_0^{t\wedge\tau_{u}^M}\Big[ 2\langle \mathcal{A}(s),  u(s)\rangle + 2\langle A(s,u(s)), X(s) \rangle - 2\langle A(s,u(s)), u(s) \rangle \nonumber\\
	& + 2 \big( \mathcal{B}(s), B(s,u(s)) \big)_{L_2} - \Vert B(s,u(s)) \Vert_{L_2}^2 \Big] ds dt .
\end{align}
Due to (\ref{220122.1716}), the limit, as $M\rightarrow\infty$, of the second term on the right hand side of the above inequality is finite.
On the other hand, by (\ref{220121.2329}) and Ito's formula we have
\begin{align}\label{220122.2057}
	& \mathbb{E} \int_0^T \psi(t) \Big[ \Vert X (t\wedge\tau_u^M) \Vert_H^2 - \Vert  x \Vert_H^2 \Big]  dt \nonumber\\
	= & \mathbb{E}  \int_0^T \psi(t) \int_0^{t\wedge\tau_{u}^M} \Big[ 2\langle \mathcal{A}(s), X(s)\rangle + \Vert \mathcal{B}(s) \Vert_{L_2}^2  \Big] ds dt .
\end{align}
Combining (\ref{220122.2055}), (\ref{220122.2056}) and (\ref{220122.2057}) together yields
\begin{align}\label{eq-8}
	& \lim_{M\rightarrow\infty}\mathbb{E} \int_0^T \psi(t) \int_0^{t\wedge\tau_{u}^M} \Big[ 2\langle \mathcal{A}(s) - A(s,u(s)), X(s) - u(s)  \rangle + \Vert \mathcal{B}(s) - B(s,u(s)) \Vert_{L_2}^2 \Big] ds dt \nonumber\\
	\leq & \lim_{M\rightarrow\infty}\liminf_{n\rightarrow\infty} \mathbb{E} \int_0^T \psi(t) \int_0^{t\wedge\tau_{u}^M} \big[ f(s) + \rho(X_n(s)) + \eta(u(s)) \big] \Vert X_n(s) - u(s) \Vert_H^2 ds  dt .
\end{align}
By the dominated convergence theorem, we can remove the limit sign on the left side of (\ref{eq-8}) to obtain
\begin{align}\label{220122.2342}
	& \mathbb{E}\int_0^T \psi(t) \int_0^{t} \Big[ 2\langle \mathcal{A}(s) - A(s,u(s)), X(s) - u(s)  \rangle + \Vert \mathcal{B}(s) - B(s,u(s)) \Vert_{L_2}^2 \Big] ds dt \nonumber\\
	\leq & C \lim_{M\rightarrow\infty}\liminf_{n\rightarrow\infty} \mathbb{E} \int_0^{T\wedge\tau_{u}^M} \big[ f(s) + \rho(X_n(s)) + \eta(u(s)) \big] \Vert X_n(s) - u(s) \Vert_H^2 ds .
\end{align}

\noindent Now take $u=X$ in the above inequality to get
\begin{align}\label{eq-9}
	& \mathbb{E}\int_0^T \psi(t) \int_0^t\Vert \mathcal{B}(s) - B(s,X(s)) \Vert_{L_2}^2 ds dt \nonumber\\
	\leq & C \lim_{M\rightarrow\infty}\liminf_{n\rightarrow\infty} \mathbb{E} \int_0^{T\wedge\tau_{X}^M} \big[ f(s) + \rho(X_n(s)) + \eta(X(s)) \big] \Vert X_n(s) - X(s) \Vert_H^2 ds .
\end{align}
Set
\begin{align}
	I := & \mathbb{E} \int_0^{T\wedge\tau_{X}^M}  f(s) \Vert X_n(s) - X(s) \Vert_H^2 ds ,\nonumber\\
	II := & \mathbb{E} \int_0^{T\wedge\tau_{X}^M}  \rho(X_n(s)) \Vert X_n(s) - X(s) \Vert_H^2 ds ,\nonumber\\
	III := & \mathbb{E} \int_0^{T\wedge\tau_{X}^M}  \eta(X(s)) \Vert X_n(s) - X(s) \Vert_H^2 ds .
\end{align}
Thus to obtain $\mathcal{B}(\cdot) = B(\cdot, X(\cdot))$, it suffices to show that
\begin{align}\label{220121.2311}
	\lim_{M\rightarrow\infty}\liminf_{n\rightarrow\infty} (I + II + III) = 0 .
\end{align}
By (\ref{220121.1753}) and Lemma \ref{220121.1807}, we have
\begin{align}\label{220210.2135}
	\lim_{M\rightarrow\infty}\lim_{n\rightarrow\infty}  I \leq \lim_{n\rightarrow\infty} \mathbb{E} \int_0^{T}  f(s) \Vert X_n(s) - X(s) \Vert_H^2 ds = 0 .
\end{align}
By (\ref{220121.1756}), the definition of $\tau_X^M$ and Lemma \ref{220121.1807}, we see that the family
\begin{align*}
  \{ \eta(X(s)) \mathbf{1}_{\{s\leq \tau_X^M\}} \Vert X_n(s) - X(s) \Vert_H^2 \}_{n=1}^{\infty}
\end{align*}
is uniformly integrable in $(s,\omega)\in [0,T]\times\Omega$. Thus by (\ref{220121.1753}), we get for any fixed $M>0$,
\begin{align*}
\lim_{n\rightarrow\infty} \mathbb{E} \int_0^{T\wedge\tau_{X}^M}  \eta(X(s)) \Vert X_n(s) - X(s) \Vert_H^2 ds = 0 .
\end{align*}
As a consequence,
\begin{align}
	\lim_{M\rightarrow\infty}\lim_{n\rightarrow\infty}  III  = 0 .
\end{align}

Next, we look at the term $II$. By (\ref{220120.1632}),
\begin{align}\label{220121.2144}
	II \leq  & C\mathbb{E} \int_0^{T\wedge\tau_X^M} (1+\Vert X_n(s)\Vert_H^{\lambda}) \Vert X_n(s) - X(s) \Vert_H^2 ds  \nonumber\\
	& + C\mathbb{E} \int_0^{T\wedge\tau_X^M} \Vert X_n(s) \Vert_V^{\theta} \Vert X_n(s)\Vert_H^{\gamma} \Vert X_n(s) - X(s) \Vert_H^2 ds \nonumber\\
	& + C\mathbb{E} \int_0^{T\wedge\tau_X^M} \Vert X_n(s) \Vert_V^{\theta} \Vert X_n(s) - X(s) \Vert_H^2 ds \nonumber\\
	=& : II_1 + II_2 + II_3.
\end{align}
Take $p$ so that
\begin{align}\label{220124.2206}
	1+ \chi < p < 1+ \frac{2L_A}{L_B} .
\end{align}
In view of (\ref{220122.2146}), we have
\begin{align}\label{220121.2146}
\begin{cases}
	\lambda +2 < p , \quad & \text{when }  \alpha\leq 2 , \\
	\lambda +2 < \alpha+p-2, \quad & \text{when } \alpha > 2 .
\end{cases}
\end{align}
If $\alpha\leq 2$, let $q=\frac{p}{\lambda+2}>1$, then Lemma \ref{220121.1807} implies that
\begin{align}\label{220121.2259}
	& C\mathbb{E}\int_0^{T} \Big[ (1+\Vert X_n(s)\Vert_H^{\lambda}) \Vert X_n(s) - X(s) \Vert_H^2 \Big]^{q} ds \nonumber\\
	\leq & C \left\{ 1+ \sup_{n\in\mathbb{N}}\mathbb{E} \Big[\sup_{s\in [0,T]} \Vert X_n(s) \Vert_H^p \Big] + \mathbb{E} \Big[\sup_{s\in[0,T]} \Vert X(s) \Vert_H^p \Big] \right\} <\infty .
\end{align}
Thus, by (\ref{220121.1753}) and the above inequality it holds that
\begin{align}\label{220121.1957}
	\lim_{M\rightarrow\infty}\lim_{n\rightarrow\infty} II_1 \leq C \lim_{n\rightarrow\infty}\mathbb{E} \int_0^{T} (1+\Vert X_n(s)\Vert_H^{\lambda}) \Vert X_n(s) - X(s) \Vert_H^2 ds   = 0 .
\end{align}
If $\alpha>2$, let $q= \frac{\alpha+p-2}{\lambda+2}>1$, then by H\"{o}lder's inequality and Lemma \ref{220121.1807}, we get
\begin{align}\label{220121.2300}
	& \mathbb{E}\int_0^T \Big[ (1+\Vert X_n(s)\Vert_H^{\lambda}) \Vert X_n(s) - X(s) \Vert_H^2 \Big]^{q} ds \nonumber\\
	\leq &  C + C \mathbb{E}\int_0^T \Vert X_n(s) \Vert_V^{\alpha} \Vert X_n(s) \Vert_H^{p-2} ds \nonumber\\
	& + C\mathbb{E} \int_0^T \Vert X_n(s) \Vert_V^{\alpha_1} \Vert X(s) \Vert_V^{\alpha_2} \Vert X_n(s) \Vert_H^{p_1} \Vert X(s) \Vert_H^{p_2} ds   \nonumber\\
	\leq & C \bigg\{ 1+\sup_{n\in\mathbb{N}}\mathbb{E} \Big[\sup_{s\in [0,T]} \Vert X_n(s) \Vert_H^p \Big] + \mathbb{E} \Big[\sup_{s\in [0,T]} \Vert X(s) \Vert_H^p \Big]  \nonumber\\
	&  + \sup_{n\in\mathbb{N}} \mathbb{E}\Big( \int_0^T \Vert X_n(s)\Vert_V^{\alpha} ds\Big)^{\frac{p}{2}} + \mathbb{E}\Big( \int_0^T \Vert X(s)\Vert_V^{\alpha} ds\Big)^{\frac{p}{2}} \bigg\} <\infty ,
\end{align}
where $\alpha_1, \alpha_2, p_1, p_2$ are nonnegative constants satisfying $\alpha_1 + \alpha_2 = \alpha$ and $p_1 + p_2 = p-2$.
Hence (\ref{220121.1957}) holds as well.
To treat the term $II_2$, we consider three cases according to the range of the parameter $\gamma$.
If
\begin{align}
	0< \gamma \leq \frac{\theta(p-2)}{\alpha} -2 ,
\end{align}
then we have
\begin{align}
	& C\mathbb{E} \int_0^{T} \Vert X_n(s) \Vert_V^{\theta} \Vert X_n(s)\Vert_H^{\gamma} \Vert X_n(s) - X(s) \Vert_H^2 ds \nonumber\\
	\leq & C \bigg\{\mathbb{E} \int_0^{T} \Big[ \Vert X_n(s) \Vert_V^{\theta} \Vert X_n(s)\Vert_H^{\gamma} \Big( \Vert X_n(s) \Vert_H^{2-\frac{\alpha-\theta}{\alpha}} + \Vert X(s) \Vert_H^{2-\frac{\alpha-\theta}{\alpha}} \Big) \Big]^{\frac{\alpha}{\theta}} ds \bigg\}^{\frac{\theta}{\alpha}} \nonumber\\
	& \times  \bigg\{\mathbb{E} \int_0^{T} \Vert X_n(s) - X(s) \Vert_H ds \bigg\}^{\frac{\alpha-\theta}{\alpha}} .
\end{align}
Similar to (\ref{220121.2300}), by H\"{o}lder's inequality and Lemma \ref{220121.1807}, we can see that in this case,
\begin{align}\label{220122.2339}
	\lim_{M\rightarrow\infty}\lim_{n\rightarrow\infty} II_2 \leq C \bigg\{ \lim_{n\rightarrow\infty}\mathbb{E} \int_0^{T} \Vert X_n(s) - X(s) \Vert_H ds \bigg\}^{\frac{\alpha-\theta}{\alpha}}  = 0 .
\end{align}
%
%
%
%
%
If
\begin{align}
	\gamma >  \frac{\theta(p-2)}{\alpha},
\end{align}
then Lemma \ref{220121.1807} implies that
\begin{align}\label{220121.2257}
	II_2 \leq & C \bigg\{ \mathbb{E} \int_0^T \Big[ \Vert X_n(s) \Vert_V^{\theta} \Vert X_n(s)\Vert_H^{\frac{\theta(p-2)}{\alpha}} \Big]^{\frac{\alpha}{\theta}} ds \bigg\}^{\frac{\theta}{\alpha}} \nonumber\\
	& \times \bigg\{ \mathbb{E} \int_0^T \Big[ \Vert X_n(s)\Vert_H^{\gamma-\frac{\theta(p-2)}{\alpha}} \Vert X_n(s) - X(s) \Vert_H^2 \Big]^{\frac{\alpha}{\alpha-\theta}} ds \bigg\}^{\frac{\alpha-\theta}{\alpha}}  \nonumber\\
	\leq & C  \bigg\{ \mathbb{E} \int_0^T \Big[ \Vert X_n(s)\Vert_H^{\gamma-\frac{\theta(p-2)}{\alpha}} \Vert X_n(s) - X(s) \Vert_H^2 \Big]^{\frac{\alpha}{\alpha-\theta}} ds \bigg\}^{\frac{\alpha-\theta}{\alpha}} .
\end{align}
In view of (\ref{220124.2206}) and (\ref{220122.2146}), we have
\begin{align}\label{220122.2313}
\begin{cases}
	2 + \gamma + \frac{2\theta}{\alpha}  < p , \quad & \text{when }  \alpha\leq 2 , \\
	4+ \gamma+ \theta - \alpha  < p, \quad & \text{when } \alpha > 2 .
\end{cases}
\end{align}
By the similar arguments  as for (\ref{220121.2259})-(\ref{220121.2300}), we can show that
\begin{align}\label{220121.2306}
	\lim_{M\rightarrow\infty}\lim_{n\rightarrow\infty} II_2 = 0 .
\end{align}
The case that
\begin{align}\label{220122.2334}
	\frac{\theta(p-2)}{\alpha} -2< \gamma \leq  \frac{\theta(p-2)}{\alpha}
\end{align}
is similar, but simpler, we omit the details. Also similar arguments lead to
\begin{align}\label{220121.2310}
	\lim_{M\rightarrow\infty}\lim_{n\rightarrow 0} II_3 = 0 .
\end{align}
Putting (\ref{220121.2144}), (\ref{220121.1957}), (\ref{220122.2339}), (\ref{220121.2306}) and (\ref{220121.2310}) together yields
\begin{align}
	\lim_{M\rightarrow\infty}\lim_{n\rightarrow\infty} II = 0 .
\end{align}
Therefore, (\ref{220121.2311}) follows and hence $\mathcal{B}(\cdot) = B(\cdot, X(\cdot))$, a.e..
\vskip 0.3cm
Taking $u=X-\varepsilon \phi e$ in (\ref{220122.2342}) for any $\varepsilon>0$, $\phi\in L^{\infty}(\Omega\times[0,T],\mathbb{\mathbb{R}})$ and $e\in V$, then dividing both sides by $\varepsilon$ and letting $\varepsilon\rightarrow 0+$, by (H1), (\ref{220120.1632}), (\ref{220121.1756}) and Lemma \ref{220121.1807}, we obtain
\begin{align}
	\mathbb{E}\int_0^T \psi(t) \int_0^t \langle \mathcal{A}(s) - A(s,X(s)), e \rangle \phi(s)  ds dt \leq 0 .
\end{align}
By the arbitrariness of $e,\phi$ and $\psi$, we conclude that $\mathcal{A}(\cdot) = A(\cdot, X(\cdot))$, a.e..
$\blacksquare$

\vskip 0.6cm
\noindent {\bf Completion of the proof of Theorem \ref{220120.1643}}.
By (\ref{220121.2329}) and Lemma \ref{220121.2330}, we know that $X$ is a probabilistically weak solution of equation (\ref{1.a}). According to   Theorem \ref{1231.0755}, the pathwise uniqueness of solutions holds. Thus the well-known Yamada-Watanabe theorem implies that there exists a unique probabilistically strong solution to equation (\ref{1.a}). The proof of the  continuity of the solution with respect to the initial value is the same as in Section \ref{220406.2023}.

\vskip 0.6cm

\section{Applications}\label{220406.2046}
\setcounter{equation}{0}

The results of in Section 2 and Section 3 are applicable to a large class of SPDE. It should be pointed out that all the examples considered in \cite{PR, LR2, LR, Liu} can be covered by our framework, including the 2D Navier-Stokes equations, porous media equations, fast-diffusion equations, $p$-Laplacian equations, Burgers equations, Allen-Cahn equations, 3D Leray-$\alpha$ model, 2D Boussinesq system, 2D magneto-hydrodynamic equations, 2D Boussinesq model for the B{\'e}nard convection, 2D magnetic B{\'e}nard equations, some shell models of turbulence (GOY, Sabra, dyadic), power law fluids, the Ladyzhenskaya model, the Kuramoto-Sivashinsky equations and the 3D tamed Navier-Stokes equations.
	In this section, we will present some examples which can not be covered in the framework previously in the literature, 		
	 but are covered by our frameworks in Section \ref{220406.2023} or Section \ref{220123.1616}.

\begin{example} [Quasilinear SPDEs] \label{220325.1642}
Let $\mathcal{O}$ be a bounded domain in $\mathbb{R}^d$ with smooth boundary $\partial\mathcal{O}$. We consider the following quasilinear partial differential equation:
\begin{align}\label{220310.1123}
	\partial_t u(t,x) = \nabla \cdot a\big(t,x,u(t,x),\nabla u(t,x)\big) - a_0\big(t,x,,u(t,x),\nabla u(t,x)\big),
\end{align}
with the zero Dirichlet boundary conditions (the case of other boundary conditions is similar),
where $u: [0,T]\times\mathcal{O}\rightarrow\mathbb{R}$,
the vector $\nabla u(t,x) = (\partial_i u(t,x))_{i=1}^d$ is the gradient of $u$ with respect to the spatial variable $x$. $a=(a_1, a_2, \cdots, a_d)$ is a vector with $a_i : [0,T]\times \mathcal{O} \times \mathbb{R}\times\mathbb{R}^{d} \rightarrow \mathbb{R}^d$ for each  $i=0,1,\cdots, n$.

We assume that $a_i$, $i=0,1,2,\cdots , d$, satisfy the following conditions: there exists a constant $\alpha>1$ if $d=1,2$ and $\alpha \geq \frac{2d}{d+2}$ if $d\geq 3$, such that
\begin{itemize}
	\item [(S1)] $a_i$ satisfies the Carathéodory conditions: for a.e. fixed $(t,x)\in [0,T]\times \mathcal{O}$, $a_i(t,x,u,z)$ is continuous in $(u,z)\in\mathbb{R}\times\mathbb{R}^{d}$, for each fixed $(u,z)\in\mathbb{R}\times\mathbb{R}^{d}$, $a_i(t,x,u,z)$ is measurable with respect to $(t,x)\in [0,T]\times\mathcal{O}$.
	\item [(S2)] There exist nonnegative constants $c_1, c_2$ and a function $f_1\in L^{\frac{\alpha}{\alpha-1}}([0,T]\times\mathcal{O}, \mathbb{R}_{+})$ such that for a.e. $(t,x)\in [0,T]\times\mathcal{O}$ and all $(u,z)\in\mathbb{R}\times\mathbb{R}^{d}$, $i=1,\cdots,d$,
		\begin{align}
			 | a_i(t,x,u,z)| \leq c_1 |z|^{\alpha -1} + c_2 |u|^{\frac{(\alpha -1)(d+2)}{d}} + f_1(t,x) .
		\end{align}
	\item [(S3)] There exist positive constants $c_3, c_4$, and a function $f_2\in L^{1}([0,T]\times\mathcal{O}, \mathbb{R}_{+})$ such that for a.e. $(t,x)\in [0,T]\times\mathcal{O}$ and all $(u,z)\in\mathbb{R}\times\mathbb{R}^{d}$,
		\begin{align}
			\sum_{i=1}^{d} a_i(t,x,u,z) z_i + a_0(t,x,u,z)u \geq c_3 |z|^{\alpha} - c_4 |u|^2 - f_2(t,x) .
		\end{align}
	\item [(S4)] For a.e. $(t,x)\in [0,T]\times\mathcal{O}$ and all $u\in\mathbb{R}$ and $z, \tilde{z}\in\mathbb{R}^{d}$ such that $z\neq \tilde{z}$,
		\begin{align}
			\sum_{i=1}^{d} [a_i(t,x,u,z) - a_i(t,x,u,\tilde{z})](z_i - \tilde{z}_i) > 0.
		\end{align}
		And for a.e. $(t,x)\in [0,T]\times\mathcal{O}$, and any $M>0$,
		\begin{align}
			\lim_{|z|\rightarrow\infty}\dfrac{\sup_{|u|\leq M} \sum_{i=1}^d a_i (t,x,u,z) z_i}{|z| + |z|^{\alpha - 1}} = \infty .
		\end{align}

	\end{itemize}

Set $H := L^2(\mathcal{O})$ and $V := W^{1,\alpha}_{0}(\mathcal{O})$, the usual Sobolev space with zero trace. By the Sobolev embedding theorem, we have the Gelfand triple
$ V  \subseteq H  \subseteq V^* $,
and the embedding $V\subseteq H$ is compact.
For $u,v\in V$, the operator $A$ is defined as follows
\begin{align}\label{220325.2322}
	\langle A(t,u), v \rangle = - & \int_{\mathcal{O}} \Big\{ \sum_{i=1}^{d} a_i\big(t,x,u(x),\nabla u(x)\big) \partial_i v(x) \nonumber\\
	& ~~~~~~ + a_0\big(t,x,,u(x),\nabla u(x)\big) v(x) \Big\} dx .
\end{align}
Recall the Gagliardo-Nirenberg inequality for $1\leq p \leq \infty$,
\begin{align}\label{220309.2230}
	\Vert u \Vert_{L^p(\mathcal{O})} \leq C \Vert \nabla u \Vert_{L^{\alpha}(\mathcal{O})}^{\delta} \Vert u \Vert_{L^2(\mathcal{O})}^{1-\delta} ,
\end{align}
where
\begin{align}
	\delta\in [0,1] \quad \text{and}\quad \frac{1}{p}=\big(\frac{1}{\alpha}-\frac{1}{d}\big)\delta + \frac{1-\delta}{2} .
\end{align}
Then it follows from (S2) and (\ref{220309.2230}) that $A$ is a measurable mapping from $[0,T]\times V$ to $V^*$, and moreover,
\begin{align}
	\Vert A(t,u) \Vert_{V^*}^{\frac{\alpha}{\alpha-1}} \leq c_1 \Vert u\Vert_{V}^{\alpha} + c \Vert u \Vert_{V}^{\alpha} \Vert u \Vert_{H}^{\frac{2\alpha}{d}} + F(t) ,
\end{align}
where
\begin{align}
	F(t) = \int_{\mathcal{O}} f_1(t,x)^{\frac{\alpha}{\alpha -1}} dx .
\end{align}
is integrable on $[0,T]$. Thus, the growth condition (H4) in Section \ref{220406.2023} is satisfied. By (S1) and (S2) it's easy to see that the hemicontinuity condition (H1) is satisfied.  (S3) and (\ref{220309.2230}) imply the coercivity condition (H3) in Section \ref{220406.2023}. By (S1), (S2) and (S4), we can show that the operator $A$ is pseudo-monotone for a.e. $t\in [0,T]$, see Theorem 10.65 and Theorem 10.63 in \cite{RR}, or Theorem 2.8 in \cite{Lions}.
Therefore, we can apply Corollary \ref{220310.1121} to obtain the existence of probabilistically weak solutions to the corresponding stochastic quasilinear partial differential equations.

\vskip 0.3cm

A typical example of (\ref{220310.1123}) is the $p$-Laplacian for $p\geq 2$,
\begin{align}
	\partial_{t}u = \nabla \cdot (|\nabla u|^{p-2} \nabla u) - c |u|^{p-2} u ,
\end{align}
where $c>0$. In this case we take $\alpha=p$, and it is easy to verify that (S1)-(S4) are satisfied.

\vskip 0.3cm

To get the uniqueness of solutions to (\ref{220310.1123}), we need to replace the assumption (S4) by the following condition:  for $\alpha\geq d$.
\begin{itemize}
	\item [(S4)$^{\prime}$] Let $0 \leq \gamma \leq \alpha(1+\frac{2}{d})-2$ and $f_3\in L^1([0,T],\mathbb{R}_{+})$. There exists a constant $c>0$ such that for a.e. $(t,x)\in [0,T]\times\mathcal{O}$ and all $u,\tilde{u}\in\mathbb{R}$ and $z, \tilde{z}\in\mathbb{R}^{d}$,
		\begin{align}
			& \sum_{i=1}^{d} [a_i(t,x,u,z) - a_i(t,x,\tilde{u},\tilde{z})](z_i - \tilde{z}_i) \nonumber\\
			 + & [a_0(t,x,u,z)-a_0(t,x,\tilde{u},\tilde{z})](u-\tilde{u}) \geq - c ( f_3(t) + |u|^{\gamma} + |\tilde{u}|^{\gamma})|u-\tilde{u}|^{2}.
		\end{align}
\end{itemize}
Under the condition (S4)$^{\prime}$, it follows from (\ref{220309.2230}) that the operator $A$ satisfies
\begin{align}
	& \langle A(t,u) - A(t,v) , u - v \rangle \nonumber\\
	\leq & \big[f_3(t) + C \Vert u\Vert_V^{\gamma\delta} \Vert u \Vert_H^{\gamma(1-\delta)} + C \Vert v\Vert_V^{\gamma\delta} \Vert v\Vert_H^{\gamma(1-\delta)} \big] \Vert u -v\Vert_H^2
\end{align}
with $\delta= \frac{\alpha d}{\alpha d + 2\alpha - 2d}$. Thus in this case, the local monotonicity  condition (H2) in Section \ref{220406.2023} is satisfied, which gives the uniqueness.


\end{example}


\begin{example}[Convection diffusion equation]\label{220325.2117}

The convection–diffusion equation describes physical phenomena where particles, energy, or other physical quantities are transferred inside a physical system due to two processes: diffusion and convection. And it has significant applications in fluid dynamics, heat transfer, and mass transfer. The convection diffusion equation is given by
\begin{align}\label{220325.1542}
\begin{cases}
	du = \nabla \cdot [a(u)\nabla u + b(u)] dt, \quad \text{on }  (0,T]\times\mathbb{T}^d ,  \\
	u(0) = u_0 ,
\end{cases}
\end{align}
where $\mathbb{T}^d$ denotes the $d$-dimensional torus, $u:[0,T]\times\mathbb{T}^d\rightarrow\mathbb{R}$, the flux function $b=(b_1,\cdots,b_d): \mathbb{R}\rightarrow\mathbb{R}^d$, the diffusion matrix $a=(a_{ij})_{i,j=1}^d: \mathbb{R}^d \rightarrow \mathcal{M}_{d\times d}$, here $\mathcal{M}_{d\times d}$ is the set of all $d\times d$-dim matrices.
We assume that $a$ and $b$ are continuous, $b$ has linear growth, $a$ is bounded and uniformly positive definite, i.e. there exists constants $\delta, C>0$ such that for any $u\in\mathbb{R}$ and $z\in\mathbb{R}^d$,
\begin{align}
 \delta |z|^2 \leq \langle a(u)z, z\rangle  \leq C |z|^2.
\end{align}

We would like to point out that under the above conditions, equation (\ref{220325.1542}) fulfills the conditions (S1) (S2) (S3) and (S4) in Example \ref{220325.1642}, but not (S4)$^{\prime}$.

In the following, we will show that equation (\ref{220325.1542}) falls into the framework in Section \ref{220406.2023}.

\vskip 0.3cm

Set $H := L^2(\mathbb{T}^d)$ and $V := W^{1,2}(\mathbb{T}^d)$. Then we have the Gelfand triple $V \subseteq H \subseteq V^*$,
and the embedding $V\subseteq H$ is compact.
For $u,v\in V$, define the operator $A$ as
\begin{align}
	\langle A(u), v \rangle = - \int_{\mathbb{T}^d} \langle a(u(x)) \nabla u(x) + b(u(x)), \nabla v(x) \rangle dx .
\end{align}
Under the above conditions on the coefficients $a$ and $b$, it is easy to see that conditions (H1), (H3) and (H4) in Section \ref{220406.2023} are satisfied, but (H2) does not hold. However, we will show  that the operator $A$ is pseudo-monotone, i.e. if
\begin{align}\label{220312.1324}
	u_n \rightharpoonup u \quad \text{weakly in } V \quad \text{and } \quad \liminf_{n\rightarrow\infty} \langle A(u_n), u_n - u\rangle \geq 0 ,
\end{align}
then for any $v\in\mathbb{V}$,
\begin{align}\label{220312.1417}
	\limsup_{n\rightarrow\infty} \langle A(u_n), u_n - v\rangle \leq \langle A(u), u-v\rangle .
\end{align}

The compact embedding of $V\subseteq H$ implies that if $u_n$ weakly converges to $u$ in $V$, then $\Vert u_n - u \Vert_H\rightarrow 0 $. Thus we can subtract a subsequence (still denoted by $\{u_n\}$) such that $u_n(x) \rightarrow u(x)$ for a.e. $x\in\mathbb{T}^d$.
Moreover, the continuity of $b$ implies that $\Vert b(u_n) - b(u)\Vert_{H}\rightarrow 0$. So
\begin{align}\label{220325.1614}
	\lim_{n\rightarrow\infty}\int_{\mathbb{T}^d} \langle b(u_n(x)), \nabla u_n(x) - \nabla u(x)\rangle dx = 0 .
\end{align}
Similarly, by the boundedness and continuity of $a$, we have
\begin{align}\label{220325.1615}
	\lim_{n\rightarrow\infty} \int_{\mathbb{T}^d} \langle a(u_n(x)) \nabla u(x), \nabla u_n(x) - \nabla u(x) \rangle dx  = 0 .
\end{align}
Combining  (\ref{220312.1324}), (\ref{220325.1614}) and (\ref{220325.1615}) together yields
\begin{align}
	- \limsup_{n\rightarrow\infty} \int_{\mathbb{T}^d} \langle a(u_n(x))(\nabla u_n(x) - \nabla u(x)), \nabla u_n(x) - \nabla u(x) \rangle dx  \geq 0.
\end{align}
Since $a$ is uniformly positive-definite, it follows from the above inequality that
\begin{align}
	\Vert u_n - u\Vert_V \rightarrow 0 .
\end{align}
Therefore, there exists a further subsequence (still denoted by $\{u_n\}$) such that $\nabla u_n(x)\rightarrow\nabla u(x)$ for a.e. $x\in\mathbb{T}^d$. Thus, for any $v\in V$,
\begin{align}
	& \lim_{n\rightarrow\infty} \int_{\mathbb{T}^d} \langle a(u_n(x))\nabla u_n(x) + b(u_n(x)), \nabla u_n(x) - \nabla v(x) \rangle dx \nonumber\\
	= & \int_{\mathbb{T}^d} \langle a(u(x))\nabla u(x) + b(u(x)), \nabla u(x) - \nabla v(x) \rangle dx .
\end{align}
Since the right side is independent of the subsequences,  the above limit holds for the whole sequence $u_n$. Hence (\ref{220312.1417}) is proved.

For the corresponding stochastic equation associated with (\ref{220325.1542}), we assume that the  diffusion coefficient is globally Lipschitz in $H$. Thus, according to Corollary \ref{220310.1121}, we obtain the existence of probabilistically weak solutions to the corresponding stochastic equation, and estimate (\ref{1226.2116}) holds. The pathwise uniqueness of the corresponding stochastic equation can be established by an argument of Yamada-Watanabe approximation under additional assumption that coefficient $a$ and $b$ are Lipschitz, see Theorem 3.1 in \cite{HZ}.

\end{example}

\begin{remark}
	The existence of stochastic convection-diffusion equations was established in \cite{HZ} under the additional assumption that  $a$ and $b$ are Lipschitz. With the approach in this paper,  Lipschitz continuity of coefficients $a$ and $b$ is not needed for the existence of solutions.
\end{remark}

\begin{example}[Cahn-Hilliard equation]

The well-known Cahn–Hilliard equations were initially introduced in \cite{CH2} to describe phase separation in a binary alloy. It is a fundamental phase field model in material science.
The classical Cahn–Hilliard equation reads:
%

\begin{align}
\begin{cases}
		\partial_t u(t)  = - \Delta^2 u + \Delta \varphi(u),   \\
	\nabla u \cdot \nu  = \nabla(\Delta u)\cdot \nu = 0 \quad \text{on } \ \partial\mathcal{O} , \\
	u(0)  = u_0,
\end{cases}
\end{align}
where $u: [0,T]\times\mathcal{O}\rightarrow\mathbb{R}$ represents a scaled concentration, $\mathcal{O}$ is a bounded domain in $\mathbb{R}^d$ with $d=1,2,3$ and with smooth boundary $\partial{O}$, $\nu$ is the outward unit normal vector on $\partial{O}$. We assume that the nonlinear term $\varphi$ satisfies the following conditions:
$\varphi\in C^1(\mathbb{R},\mathbb{R})$ and there exist constants $C\geq 0$ and $2\leq p \leq \frac{d+4}{d}$ such that for any $x,y\in\mathbb{R}$,  $\varphi^{\prime}(x) \geq -C$, $|\varphi(x)| \leq C(1+|x|^p)$ and
\[
	 |\varphi(x) - \varphi(y)| \leq C(1+|x|^{p-1} + |y|^{p-1})|x-y| .
\]

Now let $H = L^2(\mathcal{O})$ and $V= \{u\in H^2: \nabla u \cdot \nu = \nabla(\Delta u)\cdot \nu = 0 \text{  on  } \partial\mathcal{O}\}$.
Then we have the Gelfand triple $V\subseteq H \subseteq V^*$ and the embedding $V\subseteq H$ is compact. Set
\begin{align}
	A(u) = - \Delta^2 u + \Delta \varphi(u).
\end{align}
By the condition of $\varphi$ and the Gagliardo-Nirenberg inequality, conditions (H1), (H2), (H3) and (H4) in Section \ref{220406.2023} can be verified with $\alpha=2$, see Example 5.2.27 in \cite{LR2}. And condition (H2) reads as
\begin{align}\label{231028.1721}
	\langle A(u) - A(v), u-v \rangle \leq & -\frac{1}{2} \Vert u- v\Vert_V^2 + C\Big(1+ \Vert u\Vert_V^{\frac{d(p-1)}{2}} \Vert u\Vert_H^{\frac{(4-d)(p-1)}{2}} \nonumber\\
	&+ \Vert v\Vert_V^{\frac{d(p-1)}{2}} \Vert v\Vert_H^{\frac{(4-d)(p-1)}{2}} \Big) \Vert u - v\Vert_H^2 .
\end{align}
Note that $\frac{d(p-1)}{2}\leq 2 \Longleftrightarrow p\leq \frac{d+4}{d}$.
In the case of $d=1$ and $d=2$, the function $\varphi$ can be taken to be the typical example $\varphi(x) = x^3 - x$, which is the derivative of the double well potential $F(x) = \frac{1}{4}(x^2 - 1)^2$.

\vskip 0.3cm

Under the above conditions on $\varphi$, by Theorem \ref{1227.2229} and Theorem \ref{220119.1651}, we have established the well-posedness of solutions to the corresponding stochastic Cahn-Hilliard equation,
\begin{align}\label{231028.2055}
	d u(t) = \big[ - \Delta^2 u + \Delta \varphi(u) \big] dt + B(t,u) dW(t) ,
\end{align}	
where $W$ is a cylindrical Wiener process on another separable Hilbert space $U$, $B$ is Lipschitz from $H$ to $L_2(U,H)$, and the initial value $u(0)\in H$.
Moreover, in the case of $d=1$ and $p<5$ (the exponent $p$ in the typical example $\varphi(x) = x^3 - x$ is $3$), the coefficient $B$ can also depend on $\Delta u$. Because in this case, (H2)* is satisfied due to the fact that the exponent $\theta = \frac{p-1}{2}$ in (4.26) is less than $\alpha=2$. Note that (H3)* and (H4)* hold as before. Therefore, in the case of $d=1$ and $p<5$, the well-posedness of solutions to (\ref{231028.2055}) is established according to Theorem \ref{220120.1643} in Section \ref{220123.1616}, when the coefficient $B$ satisfies (H5)* and (\ref{231028.2048}) holds.
%
%
%
%
%
To the best of our knowledge, these results are not seen in literature. We refer the reader to \cite{C, LR3, AKM} and references therein.

\end{example}


\begin{example}[2D Liquid crystal model]
The elementary form of the hydrodynamics of liquid crystals is a simplified version of the Ericksen–Leslie system with Ginzburg–Landau approximation, which is established by Lin an Liu in \cite{LL}. This model in two dimensions is given by
\begin{align}\label{220323.1751}
\begin{cases}
	\partial_t u  =  \Delta u  - (u\cdot \nabla)u -\nabla p - \nabla \cdot (\nabla n \otimes \nabla n) , \\
	\nabla \cdot u = 0 , \\
	\partial_t n  = \Delta n - (u\cdot\nabla)n - \Phi(n) , \\
	 u = 0  \quad\text{and}\quad  \frac{\partial n}{\partial \nu} = 0  \quad\text{on } \partial\mathcal{O},  \\
	 u(0) = u_0, \quad n(0) = n_0 ,
\end{cases}	
\end{align}
where $\mathcal{O}$ is a bounded domain in $\mathbb{R}^2$ with smooth boundary $\partial\mathcal{O}$, $u: [0,T]\times\mathcal{O}\rightarrow \mathbb{R}^2$ is the velocity, $p: [0,T]\times\mathcal{O}\rightarrow \mathbb{R}$ is the pressure, $n: [0,T]\times\mathcal{O}\rightarrow \mathbb{R}^3$ is the director field of liquid crystal molecules, $\nu$ is the outward unit normal vector on $\partial{O}$. By the symbol $\nabla n \otimes \nabla n$ we mean a $2\times 2$ matrix with entries defined by
\[ (\nabla n \otimes \nabla n)_{i,j} = \sum_{k=1}^{3} \big(\partial_{i} n_{k} \big) \big( \partial_{j} n_{k} \big), \]
where $\partial_i$ denotes the partial derivative with respect to $x_i$ for $i=1,2$.
We assume that
$\Phi: \mathbb{R}^3 \rightarrow \mathbb{R}^3$ satisfies the following conditions:
there exists a $k$-th polynomial $\varphi: [0,\infty) \rightarrow \mathbb{R}$ for some $k\in\mathbb{N}$ such that
\[ \Phi(n) = \varphi(|n|^2)n = \Big( \sum_{i=0}^k a_i |n|^{2i} \Big) n , \]
where $a_i\in\mathbb{R}$ for $i=0,1,\cdots,k-1$ and $a_k >0$.

\vskip 0.3cm

Next we will verify that the above model falls into the framework of Section \ref{220406.2023} and Section \ref{220123.1616}.
Let $V=\{u\in H^1(\mathcal{O})^2: \nabla\cdot u = 0,\ u|_{\partial\mathcal{O}}=0\}$. Denote by $H$ the closure of $V$ under the $L^2$-norm $\Vert u\Vert_H^2 : = \int_{\mathcal{O}}|u(x)|^2 dx $.
Now set
\begin{align}
	\mathbb{H}: = H \times [ H^1(\mathcal{O})^3 ], \quad \mathbb{V}:= V\times \Big\{n\in H^2(\mathcal{O})^3 : \frac{\partial n}{\partial \nu} = 0 \Big\} ,
\end{align}
with the norm in $H$ and in $V$ denoted separately by
\[ \Vert X \Vert_{\mathbb{H}}^2 := \Vert u\Vert_H^2 + \Vert n\Vert_{H^1}^2, \quad \Vert X \Vert_{\mathbb{V}}^2 := \Vert u\Vert_V^2 + \Vert n\Vert_{H^2}^2 \]
for $X=(u,n)$.
Then we have the Gelfand triple $\mathbb{V}\subseteq \mathbb{H} \subseteq \mathbb{V}^*$ and the embedding $\mathbb{V}\subseteq \mathbb{H}$ is compact.

Note that
\begin{align}\label{220322.1001}
	 \nabla\cdot (\nabla n\otimes \nabla n) = \frac{1}{2}\nabla (|\nabla n|^2) + \nabla n \cdot \Delta n .
\end{align}
Let $P_H: L^2(\mathcal{O})^2 \rightarrow H$ be the usual Helmholtz-Leray projection.
And we set
\begin{align}
	A(X) :=
\left(
\begin{array}{c}
	 P_H [ \Delta u -  (u\cdot \nabla) u - \nabla n \cdot \Delta n ] \\
	\Delta n - (u\cdot \nabla) n - \varphi(n)
\end{array}
\right)	.
\end{align}
It is  known (see e.g. \cite{T}) that
\begin{align}\label{220322.1026}
	\Vert P_H [\Delta u - (u\cdot\nabla)u] \Vert_{V^*}^2 \leq C(1+\Vert u\Vert_H^2)\Vert u\Vert_V^2 .
\end{align}
By (\ref{220322.1001}) and (\ref{220309.2230}), we have
\begin{align}
	\Vert P_H (\nabla n \cdot \Delta n) \Vert_{V^*}^2
\leq \Vert \nabla n \Vert_{L^4(\mathcal{O})}^4 \leq \Vert  n\Vert_{H^1}^2 \Vert n \Vert_{H^2}^2 .
\end{align}
Obviously,
\begin{align}
	\Vert \Delta n \Vert_{L^2(\mathcal{O})} \leq \Vert n \Vert_{H^2} .
\end{align}
It follows from (\ref{220309.2230}) that
\begin{align}
	\Vert (u\cdot \nabla) n \Vert_{L^2(\mathcal{O})}^2 \leq \Vert u\Vert_{L^4(\mathcal{O})}^2 \Vert \nabla n\Vert_{L^4(\mathcal{O})}^2 \leq C\Vert u\Vert_H \Vert u\Vert_V \Vert n\Vert_{H^1} \Vert n\Vert_{H^2} .
\end{align}
The condition on $\Phi$ implies
\begin{align}\label{220322.1027}
	\Vert \Phi(n) \Vert_{L^2(\mathcal{O})}^2 \leq C \Vert n \Vert_{L^{4k+2}(\mathcal{O})}^{4k+2} \leq C\Vert n \Vert_{H^1}^{4k+2} .
\end{align}
Combining (\ref{220322.1026})-(\ref{220322.1027}) together, we obtain
\begin{align}
	\Vert A(X) \Vert_{\mathbb{V}^*}^2 \leq C (1+\Vert X \Vert_{\mathbb{H}}^{4k+2}) \Vert X\Vert_{\mathbb{V}^*}^2 .
\end{align}
Thus the operator $A: \mathbb{V}\rightarrow\mathbb{V}^*$ satisfies the condition (H4) in Section \ref{220406.2023}.  By the integration by parts, we have
\begin{align}\label{220323.1459}
	{}_{\mathbb{V}^*}\langle A(X), X \rangle_{\mathbb{V}} = & {}_{V^*}\langle \Delta u - (u\cdot\nabla)u - \nabla n\cdot\Delta n , u\rangle_{V} \nonumber\\
	& + {}_{L^2}\langle \Delta n - (u\cdot \nabla) n -\Phi(n), n\rangle_{H^2} \nonumber\\
	= & -\Vert u\Vert_{V}^2 - ((u\cdot\nabla )n, \Delta n )_{L^2} - \Vert \Delta n \Vert_{L^2}^2 - \Vert \nabla n\Vert_{L^2}^2 \nonumber\\
	& + ((u\cdot\nabla )n, \Delta n )_{L^2} - {}_{L^2}\langle \Phi(n), n \rangle_{H^2} \nonumber\\
	= & - (\Vert u\Vert_{V}^2 + \Vert n\Vert_{H^2}^2 ) + \Vert n\Vert_{L^2}^2 - {}_{L^2}\langle \Phi(n), n \rangle_{H^2} .
\end{align}
The last term on the right hand side of the above inequality can be estimated as follows
\begin{align}\label{220323.1500}
	& - {}_{L^2}\langle \Phi(n), n \rangle_{H^2} \nonumber\\
	= & - \big(\Phi(n), n \big)_{L^2} - \big(\nabla \Phi(n), \nabla n \big)_{L^2}  \nonumber\\
	= & - \int_{\mathcal{O}} \varphi(|n|^2)|n|^2 - \int_{\mathcal{O}}\sum_{j=1}^3\sum_{i=1}^2 \Big[ \varphi(|n|^2)\partial_i n_j + 2 \varphi^{\prime}(|n|^2)\sum_{l=1}^3 n_j n_l \partial_i n_l\Big] \partial_i n_j \nonumber\\
	= & -\int_{\mathcal{O}} \varphi(|n|^2) (|n|^2+ |\nabla n|^2) -\int_{\mathcal{O}} 2 \varphi^{\prime}(|n|^2)\mathrm{tr}\big(\nabla n \cdot( n\otimes n )\cdot (\nabla n)^{T} \big) \nonumber\\
	\leq & C\Vert n\Vert_{H^1}^2 ,
\end{align}
where we have used the fact that $\varphi(z)$ and $\varphi^{\prime}(z)z$ have lower bounds on the interval $[0,\infty)$.
Combining (\ref{220323.1459}) and (\ref{220323.1500}) together gives
\begin{align}
	{}_{\mathbb{V}^*}\langle A(X), X \rangle_{\mathbb{V}} \leq - \Vert X \Vert_{\mathbb{V}}^2 + C \Vert X \Vert_{\mathbb{H}}^2 .
\end{align}
Hence (H3) in Section \ref{220406.2023} is satisfied. For $X=(u,n)$ and $\widetilde{X}=(\widetilde{u},\widetilde{n})$ in $V$,
\begin{align}\label{220323.1631}
	& {}_{\mathbb{V}^*}\langle A(X) - A(\widetilde{X}), X - \widetilde{X}\rangle_{\mathbb{V}} \nonumber\\
	= & {}_{V^*}\langle \Delta u - \Delta \widetilde{u} , u - \widetilde{u}\rangle_{V} \nonumber\\
	& - {}_{V^*}\langle (u\cdot\nabla)u - (\widetilde{u}\cdot\nabla)\widetilde{u} , u - \widetilde{u}\rangle_{V} \nonumber\\
	& - {}_{V^*}\langle \nabla n\cdot\Delta n - \nabla \widetilde{n}\cdot\Delta \widetilde{n} , u - \widetilde{u}\rangle_{V} \nonumber\\
	& + {}_{L^2}\langle \Delta n -\Delta \widetilde{n}, n - \widetilde{n} \rangle_{H^2} \nonumber\\
	& - {}_{L^2}\langle (u\cdot \nabla) n - (\widetilde{u}\cdot \nabla) \widetilde{n}, n - \widetilde{n} \rangle_{H^2} \nonumber\\
	& - {}_{L^2}\langle \Phi(n) - \Phi(\widetilde{n}), n - \widetilde{n} \rangle_{H^2} \nonumber\\
	= & I + II + \cdots VI .
\end{align}
It is easy to see that
\begin{align}
	I & = - \Vert u - \widetilde{u} \Vert_{V}^2 , \\
	II & \leq \varepsilon \Vert u - \widetilde{u} \Vert_{V}^2 + C_{\varepsilon} \Vert u\Vert_{V}^2 \Vert u - \widetilde{u} \Vert_{H}^2 .
\end{align}
For terms $III$ and $V$, we have
\begin{align}
	III + V = & - {}_{V^*}\langle \nabla n\cdot\Delta (n - \widetilde{n}) , u - \widetilde{u}\rangle_{V} \nonumber\\
	& - {}_{V^*}\langle \nabla (n - \widetilde{n}) \cdot\Delta \widetilde{n} , u - \widetilde{u}\rangle_{V} \nonumber\\
	& - {}_{L^2}\langle ((u - \widetilde{u} ) \cdot \nabla) n , n - \widetilde{n} \rangle_{H^2} \nonumber\\
	& - {}_{L^2}\langle ( \widetilde{u} \cdot \nabla) (n - \widetilde{n}) , n - \widetilde{n} \rangle_{H^2} \nonumber\\
	= & J_1 + J_2 + J_3 + J_4.
\end{align}
By (\ref{220309.2230}) we have
\begin{align}
	|J_2|\leq & \Vert \nabla (n - \widetilde{n}) \Vert_{L^4} \Vert \Delta\widetilde{n} \Vert_{L^2} \Vert u - \widetilde{u}\Vert_{L^4} \nonumber\\
	\leq & C \Vert n - \widetilde{n}\Vert_{H^2}^{\frac{1}{2}} \Vert n - \widetilde{n}\Vert_{H^1}^{\frac{1}{2}} \Vert \widetilde{n} \Vert_{H^2}  \Vert u - \widetilde{u} \Vert_{V}^{\frac{1}{2}} \Vert u - \widetilde{u} \Vert_{H}^{\frac{1}{2}} \nonumber\\
	\leq & \varepsilon \Vert n - \widetilde{n} \Vert_{H^2} \Vert u - \widetilde{u} \Vert_{V} + C_{\varepsilon} \Vert \widetilde{n}\Vert_{H^2}^2 \Vert n - \widetilde{n} \Vert_{H^1} \Vert u - \widetilde{u} \Vert_{H} .
\end{align}
Using integration by parts and Young's inequality, we get
\begin{align}
|J_1 + J_3| = & |- {}_{L^2}\langle (u- \widetilde{u}) \cdot \nabla n, n-\widetilde{n} \rangle_{L^2}| \nonumber\\
\leq & \Vert u- \widetilde{u} \Vert_{L^4} \Vert \nabla n\Vert_{L^2} \Vert n-\widetilde{n} \Vert_{L^4} \nonumber\\
\leq & \Vert u- \widetilde{u} \Vert_{V}^{\frac{1}{2}} \Vert u- \widetilde{u} \Vert_{H}^{\frac{1}{2}} \Vert n \Vert_{H^1} \Vert n-\widetilde{n} \Vert_{H^1} \nonumber\\
\leq & \varepsilon \Vert u- \widetilde{u} \Vert_{V}^2 + C_{\varepsilon} \Vert n \Vert_{H^1}^{\frac{4}{3}} \Vert u- \widetilde{u} \Vert_{H}^{\frac{2}{3}} \Vert n-\widetilde{n} \Vert_{H^1}^{\frac{4}{3}}.
\end{align}
Similarly,
\begin{align}
	|J_4| = & | {}_{L^2}\langle ( \widetilde{u} \cdot \nabla) (n - \widetilde{n}) , \Delta ( n - \widetilde{n} ) \rangle_{L^2} | \nonumber\\
\leq &  \Vert \Delta (n - \widetilde{n} ) \Vert_{L^2} \Vert \nabla (n - \widetilde{n} )\Vert_{L^4} \Vert \widetilde{u}\Vert_{L^4} \nonumber\\
	\leq & C \Vert n - \widetilde{n} \Vert_{H^2}^{\frac{3}{2}} \Vert n - \widetilde{n}  \Vert_{H^1}^{\frac{1}{2}} \Vert \widetilde{u}\Vert_{V}^{\frac{1}{2}} \Vert \widetilde{u}\Vert_{H}^{\frac{1}{2}} \nonumber\\
	\leq & \varepsilon \Vert n - \widetilde{n} \Vert_{H^2}^{2} + C_{\varepsilon} \Vert \widetilde{u}\Vert_{V}^2 \Vert \widetilde{u}\Vert_{H}^2 \Vert n - \widetilde{n}  \Vert_{H^1}^{2} .
\end{align}
Obviously,
\begin{align}
	IV = - \Vert \Delta (n-\widetilde{n}) \Vert_{L^2}^2 - \Vert \nabla (n-\widetilde{n}) \Vert_{L^2}^2= - \Vert n-\widetilde{n} \Vert_{H^2}^2 + \Vert n-\widetilde{n} \Vert_{L^2}^2 .
\end{align}
The term $VI$ can be estimated as follows,
\begin{align}\label{220323.1632}
	|VI| \leq & \Vert n- \widetilde{n}\Vert_{H^2} \Vert \Phi(n) - \Phi(\widetilde{n})\Vert_{L^2} \nonumber\\
	\leq & \varepsilon \Vert n- \widetilde{n}\Vert_{H^2}^2 + C_{\varepsilon} \Vert \Phi(n) - \Phi(\widetilde{n})\Vert_{L^2}^2 \nonumber\\
	\leq & \varepsilon \Vert n- \widetilde{n}\Vert_{H^2}^2 + C_{\varepsilon} \int_{\mathcal{O}} (1+|n|^{2k}+|\widetilde{n}|^{2k})^2 |n-\widetilde{n}|^2 \nonumber\\
	\leq & \varepsilon \Vert n- \widetilde{n}\Vert_{H^2}^2 + C_{\varepsilon} (1+ \Vert n \Vert_{L^{4k+2}}^{4k} + \Vert \widetilde{n} \Vert_{L^{4k+2}}^{4k}) \Vert n - \widetilde{n}\Vert_{L^{4k+2}}^2  \nonumber\\
\leq & \varepsilon \Vert n- \widetilde{n}\Vert_{H^2}^2 + C_{\varepsilon}
(1+ \Vert n \Vert_{H^1}^{4k} + \Vert \widetilde{n} \Vert_{H^1}^{4k}) \Vert n - \widetilde{n}\Vert_{H^1}^2 .
\end{align}

\noindent Combining (\ref{220323.1631})-(\ref{220323.1632}) together and taking sufficiently small $\varepsilon > 0$, we obtain
\begin{align}
	& {}_{\mathbb{V}^*}\langle A(X) - A(\widetilde{X}), X - \widetilde{X}\rangle_{\mathbb{V}} \nonumber\\
	\leq & - \frac{1}{2} \Vert X - \widetilde{X} \Vert_{\mathbb{V}}^2 + C(1+ \Vert X\Vert_{\mathbb{H}}^{4k} + \Vert \widetilde{X}\Vert_{\mathbb{H}}^{4k} + \Vert \widetilde{X}\Vert_{\mathbb{V}}^{2}\Vert \widetilde{X}\Vert_{\mathbb{H}}^{2}) \Vert X - \widetilde{X}\Vert_{\mathbb{H}}^{2} .
\end{align}
Therefore, (H2) in Section \ref{220406.2023} is satisfied. The hemicontinuity condition (H1) can be easily verified by the dominated convergence theorem. Since condition (H2)* in Section \ref{220123.1616} is also satisfied, this model also falls into the framework of Section \ref{220123.1616}.

In \cite{BHR2}, the authors considered a stochastic version of system (\ref{220323.1751}) with noise in the equation of $u$ only depending on $u$, and with linear multiplicative noise only depending on $n$ in Stratonovich sense in the equation of $n$. Now applying Theorem \ref{1227.2229} and Theorem \ref{220119.1651} in Section \ref{220406.2023}, we can establish the well-posedness of the stochastic 2D liquid crystal equations driven by general multiplicative noise which can depend both on $u$ and $n$. Moreover, applying Theorem \ref{220120.1643} in Section \ref{220123.1616}, the noise can also depend on $\nabla u$ and $\Delta n$.
\end{example}

\begin{remark}
	In system (\ref{220323.1751}), if $n: \mathcal{O}\rightarrow \mathbb{R}$ and $\Phi: \mathbb{R}\rightarrow\mathbb{R}$ are  scalar functions, then the corresponding system is  the Allen-Cahn-Navier-Stokes model. This model can be viewed as a phase field model describing the motion of a mixture of two incompressible viscous fluids. We refer the readers to \cite{YFLS, LS, M} and references therein. The Allen-Cahn-Navier-Stokes model is also closely related to the magneto-hydrodynamic (MHD) equations, that is the Navier-Stokes equations coupled with the Maxwell equations. In particular in the case of dimension two and nonlinear term $\Phi(n)=0$, the corresponding system is equivalent to the MHD equations, see \cite{XZL}.  Both the Allen-Cahn-Navier-Stokes model and the MHD equations fall into the framework of Section \ref{220406.2023} and Section \ref{220123.1616} in our paper, the proof is same as above.
\end{remark}

\section{Appendix}
In this section we provide a proof of a criterion for the tightness of laws in the vector space $L^{p}([0,T], H)$. The following lemma is the  Theorem 5 of \cite{Si}.
\begin{lemma}\label{220124.2251}
	Let $1\leq p <\infty$. Let $V$, $H$ and $Y$ be Banach spaces satisfying $V \subseteq H \subseteq Y$. Suppose the embedding $V\subseteq H$ is compact. If  $\Upsilon$ is a bounded subset of  $L^{p}([0,T], V)$ satisfying
	\begin{align}
		\lim_{\delta\rightarrow 0+} \sup_{f\in\Upsilon} \int_0^{T-\delta} \Vert f(t+\delta) - f(t) \Vert_Y^p dt = 0 ,
	\end{align}
	then $\Upsilon$ is a relatively compact subset of $L^{p}([0,T], H)$.
\end{lemma}

Based on the above lemma, we can establish the following  criterion for the tightness of laws in $L^{p}([0,T], H)$.
\begin{lemma}\label{220117.1949}
		Let $1\leq p <\infty$. Let $V$, $H$ and $Y$ be Banach spaces satisfying $V \subseteq H \subseteq Y$. Suppose that the embedding $V\subseteq H$ is compact. Let $\{ X_n \}$ be a sequence of stochastic processes. If
		\begin{align}\label{220117.1935}
			\lim_{M\rightarrow\infty}\sup_{n\in\mathbb{N}}\mathbb{P}\left(  \int_0^T \Vert X_n(t) \Vert_V^p dt > M \right)  = 0 ,
		\end{align}
and for any $\epsilon>0$,
	\begin{align}\label{220117.1941}
		\lim_{\delta\rightarrow 0+} \sup_{n\in\mathbb{N}} \mathbb{P} \left( \int_0^{T-\delta} \Vert X_n (t+\delta) - X_n (t) \Vert_Y^p dt > \epsilon \right)  = 0 .
	\end{align}
	Then $\{ X_n \}$ is tight in $L^{p}([0,T], H)$.
\end{lemma}
\noindent {\bf Proof}. Take any $\varepsilon >0$. From (\ref{220117.1935}) it follows that there exists $M>0$ such that
\begin{align}\label{220117.2048}
	\sup_{n\in\mathbb{N}}\mathbb{P}\left(  \int_0^T \Vert X_n(t) \Vert_V^p dt > M \right)  \leq \frac{\varepsilon}{2}.
\end{align}
Set
\begin{align}
	K_M := \bigg\{ f\in L^p([0,T],H) : \int_0^T \Vert f(t) \Vert_V^p dt \leq M \bigg\}.
\end{align}
From (\ref{220117.1941}) it follows that for any $k\in\mathbb{N}$, there exists $\delta_k > 0$ such that
\begin{align}\label{220117.2049}
	\sup_{n\in\mathbb{N}} \mathbb{P} \left( \int_0^{T-\delta_k} \Vert X_n(t+\delta_k) - X_n(t) \Vert_Y^p dt > \frac{1}{k} \right)  \leq \frac{\varepsilon}{2^{k+1}} .
\end{align}
Set
\begin{align}
	\Gamma_k := \bigg\{ f\in L^p([0,T],H) : \int_0^{T-\delta_k} \Vert f(t+\delta_k) - f(t) \Vert_Y^p dt \leq \frac{1}{k} \bigg\} .
\end{align}
By Lemma \ref{220124.2251},
\begin{align}
	\Upsilon : = K_M \bigcap \bigcap_{k=1}^{\infty} \Gamma_k
\end{align}
is a relatively compact set in $L^p([0,T],H)$. (\ref{220117.2048}) and (\ref{220117.2049}) imply that
\begin{align}
	\sup_{n\in\mathbb{N}}\mathbb{P} \left( X_n \notin \Upsilon \right) \leq & \sup_{n\in\mathbb{N}}\mathbb{P} \left( X_n \notin K_M \right) + \sum_{k=1}^{\infty} \sup_{n\in\mathbb{N}} \mathbb{P} \left( X_n \notin \Gamma_k \right) \nonumber\\
	\leq & \frac{\varepsilon}{2} + \sum_{k=1}^{\infty} \frac{\varepsilon}{2^{k+1}} \leq \varepsilon .
\end{align}
Hence $\{ X_n \}$ is tight in $L^{p}([0,T], H)$.
$\blacksquare$

\vskip 0.6cm

\noindent {\bf Acknowledgement.} This work is partly supported by National Key R\&D
Program of China(No.2022YFA1006001) and by National Natural Science Foundation of China(No. 12131019, No. 12001516, No. 11721101),  the Fundamental Research Funds for the Central Universities (No. WK3470000031, No. WK3470000024).

\includepdf[pages=-]{corrections.pdf}

\end{document}